\documentclass[12pt]{article}
\usepackage{amsmath,amssymb,latexsym,graphics,enumerate}

\def\pa{\partial^\alpha}
\def\gaone{g_{\alpha_1}}
\def\faone{f_{\alpha_1}}
\def\faoned{f_{\alpha_1,\delta}}
\def\fgoned{f_{\gamma_1,\delta}}

\def\nd{\noindent}
\def\fa{f_\alpha}

\def\fid{f_{i,\delta}}
\def\foned{f_{1,\delta}}
\def\fzerod{f_{0,\delta}}
\def\ftwod{f_{2,\delta}}
\def\fthreed{f_{3,\delta}}
\def\fad{f_{\alpha,\delta}} 

\def\pa{\partial^\alpha}
\def\fgd{f_{\gamma_1,\delta}}

\numberwithin{equation}{section}

\begin{document}

\title{\bf Gain of Regularity for the KP-I Equation}
\author{Julie Levandosky\footnote{Mathematics Department, Framingham State College,
Framingham, MA 01701, jlevandosky@frc.mass.edu} \quad Mauricio Sep\a'{u}lveda\footnote{
Departamento de Ingenier\a'{i}a Matem\a'{a}tica, Universidad de Concepci\a'{o}n, Casilla 160-C, Concepci\a'{o}n, 
Chile. mauricio@ing-mat.udec.cl} \quad
Octavio Vera Villagr\a'{a}n\footnote{Departamento de Matem\a'{a}tica, Universidad del B\a'{i}o-B\a'{i}o 
Collao 1202, Casilla 5-C, Concepci\a'{o}n, Chile, overa@ubiobio.cl} }
\date{}
\maketitle
\begin{abstract}
\noindent In this paper we study the smoothness properties of
solutions to the KP-I equation. We show that the equation's
dispersive nature leads to a gain in regularity for the solution.
In particular, if the initial data $\phi$ possesses certain
regularity and sufficient decay as $x \rightarrow \infty$, then the
solution $u(t)$ will be smoother than $\phi$ for $0 < t \leq T$
where $T$ is the existence time of the solution.
\end{abstract}

\renewcommand{\theequation}{\thesection.\arabic{equation}}
\setcounter{equation}{0}\section{Introduction}
The KdV equation is a model for water wave propagation in shallow water with
weak dispersive and weak nonlinear effects.
In 1970, Kadomtsev \& Petviashvili \cite{KP} derived a two-dimensional analog to
the KdV equation.  Now known as the KP-I and KP-II equations, these equations
are given by
\[
u_{tx} + u_{xxxx} + u_{xx} + \epsilon u_{yy} + (uu_x)_x = 0
\]
where $\epsilon = \mp 1$.  In addition to being used as a model for the evolution of
surface waves \cite{AC}, the KP equation has also been proposed as a model for internal
waves in straits or channels of varying depth and width \cite{Sn}, \cite{DLW}.  The KP equation
has also been studied as a model for ion-acoustic wave propagation in isotropic media
\cite{PY}.  In this paper we consider smoothness properties of solutions to the
KP-I equation
\begin{eqnarray}
\label{e101}&  & (u_{t} + u_{xxx} + u_{x} + u\,u_{x})_{x} -
u_{yy} =0,\qquad (x,\,y)\in\mathbb{R}^{2},\quad t\in\mathbb{R}\\
\label{e102}&  & u(x,\,y,\,0)=\phi(x,\,y).
\end{eqnarray}

Certain results concerning the Cauchy problem for the KP-I equation
include the following.  Ukai \cite{Uk} proved local well-posedness for both the KP-I and KP-II equations
for initial data in $H^s(\mathbb R^2)$, $s \geq 3$, while Saut \cite{Sa} proved some local
existence results for generalized KP equations.  More recently, results concerning global well-posedness
for the KP-I equation have appeared.  In particular, see the works of Kenig \cite{Ke} and
Molinet, Saut, and Tzvetkov \cite{MST}.  
Here we consider the question of gain of regularity for solutions to the KP-I equation.

A number of results concerning gain of regularity for various nonlinear evolution equations
have appeared.  This paper uses the ideas of Cohen \cite{Co}, Kato \cite{Ka}, Craig and Goodman 
\cite{CG} and Craig, Kappeler, and Strauss \cite{CKS}.
Cohen considered the KdV equation, showing that ``box-shaped" initial
data $\phi \in L^2(\mathbb R^2)$ with compact support lead to a solution $u(t)$ which
is smooth for $t > 0$.  Kato generalized this result, showing that if the initial data $\phi$
are in $L^2((1+e^{\sigma x})\,dx)$, the unique solution $u(t) \in C^\infty(\mathbb R^2)$ for $t > 0$.
Kruzhkov and Faminskii \cite{KF} replaced the exponential weight function with a polynomial weight
function, quantifying the gain in regularity of the solution in terms of the decay at infinity
of the initial data.  Craig, Kappeler, and Strauss expanded on the ideas from these earlier
papers in their treatment of highly generlized KdV equations.

Other results on gain of regularity for linear and nonlinear dispersive equations include the works of
Hayashi, Nakamitsu, and Tsutsumi \cite{HNT1}, \cite{HNT2}, Hayashi and Ozawa \cite{HO}, Constantin and 
Saut \cite{CS}, Ponce \cite{Po}, 
Ginibre and Velo \cite{GV}, Kenig, Ponce and Vega \cite{KPV}, Vera \cite{thesi1}, \cite{Ve}
and Ceballos, Sepulveda and Vera \cite{CSV}. 


In studying propagation of singularities, it is natural to consider the bicharacteristics
associated with the differential operator.  For the KdV equation, it is known that the
bicharacteristics all point to the left for $t > 0$, and all singularities travel in that
direction.  Kato \cite{Ka} makes use of this uniform dispersion, choosing a nonsymmetric weight
function decaying as $x \rightarrow -\infty$ and growing as $x \rightarrow \infty$.
In \cite{CKS}, Craig, Kappeler and Strauss also make use of a unidirectional propagation of singularities
in their results on infinite smoothing properties for generalized KdV-type equations for
which $f_{u_{xxx}} \geq c > 0$.

For the two-dimensional case, Levandosky \cite{Le1} proves smoothing properties for the
KP-II equation.  This result makes use of the fact that the bicharacteristics all point into one
half-plane.  Subsequently, in \cite{Le2}, Levandosky considers generalized KdV-type equations in two-dimensions,
proving that if all bicharacteristics point into one half-plane, an infinite gain in regularity
will occur, assuming sufficient decay at infinity of the initial data.

In this paper, we address the question regarding gain in regularity for the KP-I
equation.  Unlike the KP-II equation, the bicharacteristics for the KP-I equation are not
restricted to a half-plane but span all of $\mathbb R^2$.  As a result, singularities may
travel in all of $\mathbb R^2$.  However, here we prove that if the initial data decays
sufficiently as $x \rightarrow \infty$, then we will gain a finite number of derivatives in
$x$ (as well as mixed derivatives).  In order to
state a special case of our gain in regularity theorem, we first introduce
certain function spaces we will be using.
\\

\nd{\it Definition.} We define
\begin{eqnarray}
\label{e106}X^{0}(\mathbb{R}^{2})=
\left\{u:\;u,\;\xi^{3}\widehat{u},\;\frac{\eta^{2}}{\xi}\,\widehat{u}\in
L^{2}(\mathbb{R}^{2})\right\}
\end{eqnarray}
equipped with the natural norm. On the space
\begin{eqnarray}
\label{e107}\widetilde{X}^{0}(\mathbb{R}^{2})=
\left\{u:\;\frac{1}{\xi}\,\widehat{u}(\xi,\,\eta)\in
L^{2}(\mathbb{R}^{2})\right\}
\end{eqnarray}
we define the operator $\partial_{x}^{-1}$ by
$\widehat{\partial_{x}^{-1}u}\equiv\frac{1}{i\,\xi}\,\widehat{u}.$
Therefore, in particular, we can write the norm of
$X^{0}(\mathbb{R}^{2})$ as
\begin{eqnarray}
\label{e108}||u||_{X^{0}(\mathbb{R}^{2})}^{2}=\int_{\mathbb{R}^{2}}[\,u^{2}
+ u_{xxx}^{2} + (\partial_{x}^{-1}u_{yy})^{2}\,]\,dx\,dy<+\infty
\end{eqnarray}
On this space of functions $X^{0}(\mathbb{R}^{2}),$ it makes sense
to rewrite \eqref{e101}-\eqref{e102} as
\begin{eqnarray}
\label{e109}&  & u_{t} + u_{xxx} + u_{x} + u\,u_{x} -
\partial_{x}^{-1}u_{yy} =0,\qquad (x,\,y)\in\mathbb{R}^{2},\quad t\in\mathbb{R}\\
\label{e110}&  & u(x,\,y,\,0)=\phi(x,\,y)
\end{eqnarray}
and consider weak solutions $u\in X^{0}(\mathbb{R}^{2}).$\\
\\
{\it Definition.} Let N be a positive integer. We define the space
of functions $X^{N}(\mathbb{R}^{2})$ as follows
\begin{eqnarray}
\label{e111}X^{N}=\left\{u:\;u\in L^{2}(\mathbb{R}^{2}),\;{\cal
F}^{-1}(\xi^{3}\,\widehat{u})\in H^{N}(\mathbb{R}^{2}),\,{\cal
F}^{-1}\left(\frac{\eta^{2}}{\xi}\,\widehat{u}\right)\in
H^{N}(\mathbb{R}^{2})\right\}
\end{eqnarray}
equipped with the norm
\begin{eqnarray}
\label{e112}||u||_{X^{N}(\mathbb{R}^{2})}^{2} =
\int_{\mathbb{R}^{2}}\left(u^{2} + \sum_{|\alpha|\leq
N}[\,(\pa u_{xxx})^{2} +
(\pa \partial_{x}^{-1}u_{yy})^{2}\,]\right)\,dx\,dy<+\,\infty
\end{eqnarray}
where $\alpha=(\alpha_1,\,\alpha_2)\in\mathbb{Z}^{+}\times\mathbb{Z}^{+}$
and $|\alpha|=\alpha_1 + \alpha_2.$\\
\\

\nd{\bf Gain of Regularity Theorem.} {\it Let $u$ be a solution of \eqref{e109}-\eqref{e110}
in $\mathbb R^2 \times [0,T]$ such that $u \in L^\infty([0,T]; X^1(\mathbb R^2))$ and
\begin{equation}
\label{weighted-bound-u}
\sup_{0 \leq t \leq T} \int_{\mathbb R^2} [u^2 + (\partial_y^L u)^2 + (1+x_+)^L (\partial_x^L
u)^2] \,dx \,dy < + \infty
\end{equation}
for some integer $L \geq 2$.  Then
\[
\sup_{0 \leq t \leq T} \int_{\mathbb R^2} t^{|\alpha|-L}(x_+^{2L-|\alpha|-\alpha_2} + e^{\sigma x_-})
(\pa u)^2 + \int_0^T \int_{\mathbb R^2} t^{|\alpha|-L}(x_+^{2L-|\alpha|-\alpha_2-1} + e^{\sigma x_-})
(\pa u_x)^2 < \infty,
\]
for $L+1 \leq |\alpha| \leq 2L-1$, $2L-|\alpha|-\alpha_2 \geq 1$, $\sigma > 0$ arbitrary.}
\\

\nd{\bf Remarks.}
\begin{enumerate}[(1)]
\item If we consider $|\alpha| = 2L-1$ above, then $\alpha = (2L-1,0)$ in which case the
result states that
\[
\sup_{0 \leq t \leq T} \int_{\mathbb R^2} t^{L-1} (x_+ + e^{\sigma x_-})(\partial_x^{2L-1}u)^2
+ \int_0^T \int_{\mathbb R^2} t^{L-1} (1 + e^{\sigma x_-})(\partial_x^{2L}u)^2 < \infty.
\]
In particular, this result shows a {\it gain} in $L$ derivatives in $x$.

\item While we gain $x$ derivatives and mixed derivatives, we do not gain pure $y$ derivatives.
However, we do not require any weighted estimates on $\partial_y^L u$.
In addition, we do not require any weighted estimates on $u$.  The results on the KP-II
equation include gains in pure $y$ derivatives, but also require weighted estimates on
$\partial_y^L u$.

\item The assumptions on $u$ are reasonable and shown to hold in section 6.
\end{enumerate}

The main idea of the proof is the following.  We use an inductive argument where on each level
$|\alpha|$, we apply the operator $\pa = \partial_x^{\alpha_1}
\partial_y^{\alpha_2}$ to \eqref{e109},
multiply the differentiated equation by $2f_\alpha\pa u$ where $f_\alpha$ is our weight function,
to be specified later, and integrate over $\mathbb R^2$.  Doing so, we arrive at the following
inequality
\begin{equation}
\begin{split}
\label{main-eq}
& \partial_t \int f (\pa u)^2 + 3 \int f_x (\pa u_x)^2 \leq \int f_x (\pa \partial_x^{-1} u_y)^2 \\
& \qquad + \int [f_t + f_{xxx} + f_x] (\pa u)^2 + \left|2 \int f (\pa u)\pa (uu_x)\right|
\end{split}
\end{equation}
where $\int = \int_{\mathbb R^2} dx dy$.  Assuming $f_x > 0$, the
second term on the left-hand side has a positive sign, thus allowing
us to prove a gain in regularity.  We notice that the first term on
the right-hand side is of order $|\alpha|$.  By choosing appropriate
weight functions for each $\alpha$, we have a bound on that term
from the previous step of the induction.  After proving estimates
involving the nonlinear term on the right-hand side of the equation,
we apply Gronwall's inequality to prove the bounds on the terms on
the left-hand side of the equation.

The plan of the paper is the following.  In section 2 we show the derivation of
\eqref{main-eq}.  In sections 3 and 4 we prove an existence result showing that for
initial data $\phi \in X^N(\mathbb R^2)$ there exists a smooth solution $u \in L^\infty([0,T];
X^N(\mathbb R^2))$ for a time $T$ depending only on $||\phi||_{X^0}$.  In section 5 we prove
estimates for the terms on the right-hand side of \eqref{main-eq}.  In section 6 we prove
a priori estimates showing that the solution $u$ found in section 4 also satisfies
\eqref{weighted-bound-u} for the same time $T$ as long as
\[
\int (\phi^2 + (\partial_y^L \phi)^2 + (1+x_+)^L(\partial_x^L \phi)^2) < \infty.
\]
Once we have found the solution $u$ in the appropriate weighted space as well as bounds for
terms on the right-hand side of \eqref{main-eq}, in section 7, we can state and prove
our main gain in regularity result.  This proof uses an inductive argument along with the
main estimates proven in section 5.
\\

\nd{\bf Choice of weight function. } We will be using non-symmetric weight functions.  In particular,
we will be using weight functions $f(x,t) \in C^\infty$ which behave roughly like powers of
$x$ for $x > 1$ and decay exponentially for $x < -1$.  We define our weight classes as follows.
\\

\nd{\it Definition.} A function $f = f(x,\,t)$ belongs to the weight
class $W_{\sigma \;i\;k}$ if it is a positive $C^{\infty }$
function on $\mathbb{R}\times [0,\,T]$ and there are
constants $c_{j},\,0\leq j\leq 5$ such that
\begin{eqnarray}
\label{e103}&0<c_{1}\leq t^{-\,k}\,e^{-\,\sigma \,x}\,f(x,\,t)
\leq c_{2}\qquad \forall
\;x<-1,\quad 0<t<T.&\\
\label{e104}&0<c_{3}\leq t^{-\,k}\,x^{-\,i}\,f(x,\,t)\leq
c_{4}\qquad \forall
\;x>1,\quad 0<t<T.&\\
\label{e105}&\left(t\mid \partial_{t}f\mid + \mid
\partial_{x}^{r}f\mid \right)/f\leq c_{5}\qquad \forall
\;(x,\,t)\in \mathbb{R}\times [0,\,T],\quad \forall \;r\in
\mathbb{N}.&
\end{eqnarray}
Thus $f$ looks like $t^{k}$ as $t\rightarrow 0,$ like $x^{i}$ as
$x\rightarrow +\infty$ and like $e^{\sigma\,x}$ as $x\rightarrow
-\infty.$ \\
\\
Before proceeding, we introduce some other function spaces we will be using.
\\

\nd{\it Definition.} Let $N$ be a positive integer. Let
$\widetilde H^{N}_x(W_{\sigma\;i\;k})$ be the space of functions
\begin{eqnarray}
\label{e113}\widetilde H^{N}_x(W_{\sigma \;i\;k})= \left\{v\colon
\mathbb{R}^{2}\rightarrow \mathbb{R}\,:\,\;\,||v||_{\widetilde
H^{N}_x(W_{\sigma \;i\;k})}^{2}= \int \sum_{|\alpha| \leq N} [(\pa
v)^2 + f|\partial_x^{N} v|^{2}] <+\,\infty \,\right\}
\end{eqnarray}
with $f\in W_{\sigma \;i\;k}$ fixed.\\
\\
{\bf Remarks.}
\begin{enumerate}[(1)]
\item We note that although the norm above depends on
$f,$ all choices of $f$ in this class lead to equivalent
norms.
\item The usual Sobolev space is $\,H^{N}(\mathbb{R}^{2})$ without a weight.
\end{enumerate}

\nd{\it Definition.} For fixed $f\in W_{\sigma \;i\;k}$ define the space
($N$ be a positive integer)
\begin{eqnarray}
\lefteqn{L^{2}([0,\,T]:\,\widetilde H_x^{N}(W_{\sigma \;i\;k}))}\nonumber \\
\label{e114}& = &
  \left\{v(x,\,y,\,t):\;
  ||v||_{L^{2}([0,\,T]:\,\widetilde H_x^{N}(W_{\sigma \;i\;k}))}^{2}= \int_{0}^{T}||v(\,\cdot
\,,\,\cdot\,,\,t)||_{\widetilde H_x^{N}(W_{\sigma \;i\;k})}^{2}dt<+\,\infty
\,\right\}
\end{eqnarray}
\begin{eqnarray}
\lefteqn{L^{\infty }([0,\,T]:\,\widetilde H_x^{N}(W_{\sigma \;i\;k})) }
\nonumber \\
\label{e115}& = & \left\{\,v(x,\,y,\,t):\;||v||_{L^{\infty
}([0,\,T]:\,\widetilde H_x^{N}(W_{\sigma \;i\;k}))}= \sup_{t\in
[0,\,T]}\,||v(\,\cdot\,,\,\cdot \,,\,t)||_{\widetilde H_x^{N}(W_{\sigma
\;i\;k})}<+\,\infty \,\right\}.
\end{eqnarray}
For simplicity, let
\begin{eqnarray}
\label{e118}{\cal Z}_{L} = X^{1}(\mathbb{R}^{2})\bigcap
\widetilde{H}_x^{L}(W_{0\;L\;0}).
\end{eqnarray}
With this notation, ${\cal Z}_{L}$ consists of those functions $u$
such that
\begin{eqnarray}
\label{e119}||u||_{{\cal
Z}_{L}}^{2}=\int_{\mathbb{R}^{2}}[u^2 +
u_{xxxx}^{2} + (\partial_{x}^{-1}u_{yy})^{2} + u_{yy}^2 + \sum_{|\alpha| \leq L}
(\pa u)^2 + f(\partial_x^L u)^2 \,]\,dx\,dy
\end{eqnarray}
for some $f\in W_{0\;L\;0}.$\\
\\
We now state a lemma describing one of the types of bounds we will be using for our
a priori estimates.
\\

\nd{\bf Lemma 1.1.} {\it For $p,\,q>1,$ $\frac{1}{p} + \frac{1}{q}<1,$
$u\in L^{2}(\mathbb{R}^{2}),$}
\begin{eqnarray}
\label{e120}||u||_{L^{\infty}(\mathbb{R}^{2})}\leq
c\left(\int_{\mathbb{R}^{2}}[\,1 + |\xi|^{p} +
|\eta|^{q}\,]\,|\widehat{u}|^{2}\,d\xi\,d\eta\right)^{1/2}.
\end{eqnarray}
{\it Proof.} The proof follows from writing $u$ in terms of its
inverse Fourier transform and using the fact that
\begin{eqnarray*}
\int_{\mathbb{R}^{2}}\frac{1}{1 + |\xi|^{p} +
|\eta|^{q}}\,d\xi\,d\eta <+\infty
\end{eqnarray*}
for $p,\,q$ satisfying our hypothesis. \hfill $\square$ \\
\\
In particular, we have:
\begin{eqnarray}
\label{e121}||u||_{L^{\infty}(\mathbb{R}^{2})}\leq
c\left(\int_{\mathbb{R}^{2}}[\,u^{2} + u_{xx}^{2} +
u_{y}^{2}\,]\,dx\,dy\right)^{1/2}.
\end{eqnarray}

\renewcommand{\theequation}{\thesection.\arabic{equation}}
\setcounter{equation}{0}\section{Main Equality} We consider the
KP-I equation
\begin{eqnarray}
\label{e201}&  & u_{t} + u_{xxx} + u_{x} + u\,u_{x} -
\partial_{x}^{-1}u_{yy} =0,\qquad (x,\,y)\in\mathbb{R}^{2},\quad t\in\mathbb{R}\\
\label{e202}&  & u(x,\,y,\,0)=\phi(x,\,y).
\end{eqnarray}
{\bf Lemma 2.1.} {\it Let $u$ be a solution of
\eqref{e201}-\eqref{e202} with enough Sobolev regularity and with
sufficient decay at infinity.  Let $f=f(x,t)$.  Then}
\begin{eqnarray}
\lefteqn{\partial_{t}\int_{\mathbb{R}^{2}}f
(\pa u)^{2}\,dx\,dy +
\int_{\mathbb{R}^{2}}g\,(\pa u_x)^{2}\,dx\,dy}\nonumber \\
\label{e203}&  & +\int_{\mathbb{R}^{2}}\theta\,
(\pa u)^{2}\,dx\,dy +
\int_{\mathbb{R}^{2}}\theta_{1}\, (\pa \partial_x^{-1} u_y)^{2}\,dx\,dy +
\int_{\mathbb{R}^{2}}R_{\alpha}\,dx\,dy = 0
\end{eqnarray}
{\it such that}
\begin{eqnarray*}
g & = & 3\,f_{x}\\
\theta & = & -\;[f_{t} + f_{xxx} + f_{x}\,]\\
\theta_{1} & = & -\;f_{x}\\
R_{\alpha} & = & 2\sum_{n=0}^{\alpha_1}\sum_{m=0}^{\alpha_2}
{\alpha_1\choose n}{\alpha_2\choose m}f\,(\pa u)\,(\partial_{x}^{n}\partial_{y}^{m}u)
\,(\partial_{x}^{\alpha_1 + 1 - n}\partial_{y}^{\alpha_2 - m}u).
\end{eqnarray*}
{\it Proof.} Applying the operator
$\pa$ to \eqref{e201},
we have
\begin{eqnarray*}
&  & \pa u_{t} + \pa u_{xxx} + \pa u_x + \pa (u\,u_x) -
\pa \partial_x^{-1} u_{yy} =0.
\end{eqnarray*}
Multiplying by $2\,f\,
\pa u$ and integrating
over $\mathbb{R}^{2},$ we have
\begin{equation}
\begin{split}
\label{e204}
& 2\int f\,(\pa u)\,(\pa u)_{t} +
2\int f\,(\pa u)\,(\pa u_{xxx})
+ 2\int f\,(\pa u)\,(\pa u_x) \\
& \qquad + 2\int f\,
(\pa u)\,\pa (u\,u_{x}) -\;2\int f\,
(\pa u)\,(\pa \partial_{x}^{-1} u_{yy}) =0.
\end{split}
\end{equation}
Each term in \eqref{e204} is calculated separately integrating by
parts
\begin{eqnarray*}
2\int f\,
(\pa u)\,(\pa u)_{t} = \partial_{t}\int f\,(\pa u)^{2} -
\int f_{t}\,(\pa u)^{2}.
\end{eqnarray*}
\begin{eqnarray*}
2\int f\,(\pa u)\,(\pa u_{xxx})
=3\int f_{x}\,(\pa u_x)^{2} -
\int f_{xxx}\,(\pa u)^{2}.
\end{eqnarray*}
\begin{eqnarray*}
2\int f\,(\pa u)\,(\pa u_x) =
-\int f_{x}\,(\pa u)^{2}
\end{eqnarray*}
\begin{eqnarray*}
-\;2\int f\,
(\pa u)\,(\pa \partial_{x}^{-1}u_{yy}) =-\int f_{x}\,
(\pa \partial_{x}^{- 1}u_y)^{2}.
\end{eqnarray*}
\begin{eqnarray*}
2\int f\,
(\pa u)\,\pa (u\,u_{x}) = \;2\sum_{n=0}^{\alpha_1}\sum_{m=0}^{\alpha_2}{\alpha_1\choose
n}{\alpha_2\choose
m}\int f\,(\pa u)\,(\partial_{x}^{n}\partial_{y}^{m}u)\,(\partial_{x}^{\alpha_1 +
1 - n}\partial_{y}^{\alpha_2 - m}u).
\end{eqnarray*}
Replacing in \eqref{e204} we obtain
\begin{equation*}
\begin{split}
& \partial_{t}\int f\,(\pa u)^{2} +
3\int f_{x}\,(\pa u_x)^{2} \\
& - \int [f_{t} + f_{xxx} + f_{x}]\,(\pa u)^{2} -
\int f_{x}\,(\pa \partial_{x}^{-
1}u_y)^{2} \\
& +\;2\sum_{n=0}^{\alpha_1}\sum_{m=0}^{\alpha_2}{\alpha_1\choose
n}{\alpha_2\choose m}\int f\,(\pa u)\,(\partial_{x}^{n}
\partial_{y}^{m}u)\,(\partial_{x}^{\alpha_1 + 1
- n}\partial_{y}^{\alpha_2 - m}u) = 0.
\end{split}
\end{equation*}
Therefore, we obtain the {\bf Main Equality},
\begin{equation*}
\begin{split}
& \partial_{t}\int f
(\pa u)^{2} +
3\int f_{x}\,(\pa u_x)^{2} + \int \theta\,
(\pa u)^{2} +
\int \theta_{1}\, (\pa \partial_{x}^{-
1}u_y)^{2} +
\int R_{\alpha} = 0
\end{split}
\end{equation*}
such that
\begin{eqnarray*}
\theta & = & -\;[f_{t} + f_{xxx} + f_{x}\,]\\
\theta_{1} & = & -\;f_{x}\\
R_{\alpha} & = & 2\sum_{n=1}^{\alpha_1}\sum_{m=1}^{\alpha_2}{\alpha_1\choose n}{\alpha_2
\choose m}f\,(\pa u)\,(\partial_{x}^{n}\partial_{y}^{m}u)\,(\partial_{x}^{\alpha_1 +
1 - n}\partial_{y}^{\alpha_2 - m}u).
\end{eqnarray*}


\renewcommand{\theequation}{\thesection.\arabic{equation}}
\setcounter{equation}{0}\section{An a priori estimate} In section four we
prove a basic local-in-time existence theorem for \eqref{e201}-\eqref{e202}.
The proof relies on approximating \eqref{e201} by a sequence of linear equations.
In this section, we prove an existence theorem for linear equations as well as an
a priori estimate on those solutions which will be necessary for our main existence theorem
in the next section.

We begin by approximating \eqref{e201} by the linear equation
\begin{eqnarray}
\label{e305}&  & u_{t}^{(n)} + u_{xxx}^{(n)} + u_{x}^{(n)} + u^{(n
- 1)}\,u_{x}^{(n)} -
\partial_{x}^{-1}u_{yy}^{(n)} =0
\end{eqnarray}
where the initial condition is given by
$u^{(n)}(x,\,y,\,0)=\phi(x,\,y)$ and the first approximation is
given by $u^{(0)}(x,\,y,\,t)=\phi(x,\,y).$ The linear equation
which is to be solved at each iteration is of the form
\begin{eqnarray}
\label{e306}&  & u_{t} + u_{xxx} + u_{x} + b\,u_{x} -
\partial_{x}^{-1}u_{yy} =0.
\end{eqnarray}
where $b$ is a smooth bounded coefficient.  Below we show that this equation
can be solved in any interval of time in which the coefficient is defined. \\

\nd{\bf Lemma 3.1} (Existence for linear equation). {\it Given initial
data $\phi\in H^{\infty}(\mathbb{R}^{2})=\bigcap_{N\geq
0}H^{N}(\mathbb{R}^{2})$ and
$\partial_{x}^{-1}\phi_{yy}\in\bigcap_{N\geq
0}H^{N}(\mathbb{R}^{2})$ there exists a unique solution of
\eqref{e306}. The solution is defined in any time interval in
which the coefficients are defined.}\\
\\
{\it Proof.} Let $T>0$ be arbitrary and $M>0$ be a constant. Let
\begin{eqnarray*}
{\cal L} = \partial_{t} + \partial_{x}^{3} + \partial_{x} +
b\,\partial_{x} - \partial_{x}^{-1}\partial_{y}^{2}
\end{eqnarray*}
be defined on those functions $u\in X^{0}(\mathbb{R}^{2}).$ Recall
that $u\in X^{0}(\mathbb{R}^{2})$ means
$u,\,u_{xxx},\,\frac{\eta^{2}}{\xi}\widehat{u}\in
L^{2}(\mathbb{R}^{2}).$ We consider the bilinear form ${\cal
B}:{\cal D}\times {\cal D}\rightarrow \mathbb{R},$
\begin{eqnarray*}
{\cal B}(u,\,v) =
\left<u,\,v\right>=\int_{0}^{T}\int_{\mathbb{R}^{2}}e^{-M\,t}\,u\,v\,dx\,dy\,dt
\end{eqnarray*}
where ${\cal D}=\{u\in
C([0,\,T]:\,L^{2}(\mathbb{R}^{2})):\,u(x,\,y,\,0)=0\}.$ By
integration by parts, we see that
\begin{eqnarray*}
\int_{\mathbb{R}^{2}}{\cal L}u\cdot u\,dx\,dy & = &
\frac{1}{2}\,\partial_{t}\int_{\mathbb{R}^{2}}u^{2}\,dx\,dy -
\frac{1}{2}\int_{\mathbb{R}^{2}}b_{x}\,u^{2}\,dx\,dy \\
& \geq & \frac 12 \partial_t \int_{\mathbb R^2} u^2 \,dx\,dy - \frac 12 \int_{\mathbb R^2}
c u^2 \,dx\,dy
\end{eqnarray*}
We multiply by $e^{-M\,t}$ and integrate in time to obtain for
$u\in C([0,\,T]:\,X^{0}(\mathbb{R}^{2}))$ with $u(x,\,y,\,0)=0$
\begin{eqnarray}
\label{e307}\left<{\cal L}u,\,u\right> \geq
e^{-M\,t}\int_{\mathbb{R}^{2}}u^{2}\,dx\,dy + (M -
c)\int_{\mathbb{R}^{2}}e^{-M\,t}\,u^{2}\,dx\,dy.
\end{eqnarray}
Thus, $\left<{\cal L}u,\,u\right>\geq \left<u,\,u\right>$ provided $M$ is chosen large
enough. Similarly, $\left<{\cal L}^{*}v,\,v\right>\geq \left<v,\,v\right>$ for all $v\in
C([0,\,T]:\,X^{0}(\mathbb{R}^{2}))$ such that $v(x,\,y,\,T)=0$ where
${\cal L}^{*}$ denotes the formal adjoint of ${\cal L}.$ Therefore,
$\left<{\cal L}^{*}v,\,{\cal L}^{*}u\right>$ is an inner product on ${\cal
D}^{*}=\{v\in
C([0,\,T]:\,X^{0}(\mathbb{R}^{2})):\,v(x,\,y,\,T)=0\}.$ Denote by
$Y$ the completion of ${\cal D}^{*}$ with respect to this inner
product. By the Riesz representation theorem, there exists a unique
solution $V\in Y,$ such that for any $v\in {\cal D}^{*},$ $\left<{\cal
L}^{*}V,\,{\cal L}^{*}v\right>=(\phi,\,v(x,\,y,\,0))$ where we used the
fact that $(\phi,\,v(x,\,y,\,0))$ is a bounded linear functional on
${\cal D}^{*}.$ Then $w={\cal L}^{*}V$ is a weak solution of ${\cal
L}w=0,$ $w_{0}=\phi$ with $w\in
L^{2}(\mathbb{R}^{2}\times [0,\,T]).$

{\it Remark.} To obtain higher regularity of the solution, we
repeat the proof with higher derivatives included in the inner
product. It is a standard approximation procedure to obtain a
result for general initial data. \hfill $\square$ \\
\\

Next, we need to introduce a new function space.  Let
\begin{eqnarray}
\label{e303}{Z}_{T}^{N}=\{u:\;u\in L^{\infty}([0,\,T]:\,H^{(N +
3,\,N + 2)}(\mathbb{R}^{2})),\;u_{t}\in
L^{\infty}([0,\,T]:\,H^{N}(\mathbb{R}^{2}))\}
\end{eqnarray}
where $H^{(\alpha_1,\,\alpha_2)}(\mathbb{R}^{2})=\{u:\;u,\,
\partial_{x}^{\alpha_1}u,\,\partial_{y}^{\alpha_2}u\in
L^{2}(\mathbb{R}^{2})\}$ with the accompanying norm
\begin{eqnarray}
\label{e304}||u||_{{Z}_{T}^{N}}^{2}=\sup_{t\in
[0,\,T]}\int_{\mathbb{R}^{2}}\left(u^{2} +
\sum_{|j|=N}[\,(\partial^{j}u_{xxx})^{2} +
(\partial^{j}u_{yy})^{2}\,]\right)+
\int_{\mathbb{R}^{2}}\left(u_{t}^{2} +
\sum_{|j|=N}(\partial^{j}u_{t})^{2}\right).
\end{eqnarray}

\nd Using this function space and the linearized equation \eqref{e305}, we consider the mapping
$\Pi: Z_{T}^{N} \rightarrow Z_{T}^{N}$ such that $u^{(n)} = \Pi(u^{(n-1)})$ and our
first approximation is given by $u^{(0)}(x,y,t) = \phi(x,y)$.
In Lemma 3.2 below, we show an a priori estimate which will be used on our sequence of
solutions $\{u^{(n)}\}$ in our main existence theorem in section four.\\
\\
{\bf Lemma 3.2.} {\it Let $v,\,w$ be a pair of functions in ${Z}_{t}^{N}$ for all $N$ and all $t\geq 0,$ such that $v,\,w$ are
solutions to}
\begin{eqnarray}
\label{e308}&  & v_{t} + v_{xxx} + v_{x} + w\,v_{x} -
\partial_{x}^{-1}v_{yy} =0.
\end{eqnarray}
{\it Then for all $N\geq 0,$ the following inequality holds:
\begin{eqnarray}
\label{e309}||v||_{{Z}_{t}^{N}}^{2} \leq
||v(\,\cdot\,,\,\cdot\,,\,0)||_{H^{(N + 3,\,N +
2)}(\mathbb{R}^{2})}^{2} +
||v_{t}(\,\cdot\,,\,\cdot\,,\,0)||_{H^{N}(\mathbb{R}^{2})}^{2} +
c\,t\,||w||_{{Z}_{t}^{N}}\;||v||_{{Z}_{t}^{N}}^{2}
\end{eqnarray}
for all $t\geq 0.$}\\
\\
{\it Proof.} We will show that for each $j,$ $|j|\geq 0$ and
$0\leq \widetilde{t}\leq t,$
\begin{align*}
\partial_{t}\int[\,(\partial^{j}v(\,\cdot\,,\,\cdot\,,\,\widetilde{t}))^{2}
+
& (\partial^{j}\partial_{x}^{3}v(\,\cdot\,,\,\cdot\,,\,\widetilde{t}))^{2}
+
(\partial^{j}\partial_{y}^{2}v(\,\cdot\,,\,\cdot\,,\,\widetilde{t}))^{2}
+
(\partial^{j}v_{t}(\,\cdot\,,\,\cdot\,,\,\widetilde{t}))^{2}\,] \\
& \leq
c\,||w||_{{Z}_{t}^{|j|}}\,||v||_{{Z}_{t}^{|j|}}^{2}.
\end{align*}
We begin by taking $j$ derivatives of \eqref{e308}. We have
\begin{eqnarray}
\label{e310}&  & \partial^{j}v_{t} + \partial^{j}v_{xxx} +
\partial^{j}v_{x} + \partial^{j}(w\,v_{x}) -
\partial^{j}\partial_{x}^{-1}v_{yy} =0.
\end{eqnarray}
Multiply \eqref{e310} by $2\,\partial^{j}v$ and integrate over
$\mathbb{R}^{2}.$ Hence
\begin{eqnarray*}
\partial_{t}\int
(\partial^{j}v(\,\cdot\,,\,\cdot\,,\,\widetilde{t}))^{2}  &
\leq &
c\,\left|\int\partial^{j}(w\,v_{x})\,(\partial^{j}v) \right|\\
& \leq &\left|\int[\,(\partial^{j}w)\,v_{x}\;+\,
\ldots\,+\;w\,(\partial^{j}v_{x})\,]\,(\partial^{j}v) \right|.
\end{eqnarray*}
The remainder terms can be bounded as follows:
\begin{eqnarray*}
\left|\int(\partial^{j}w)\,v_{x}\,(\partial^{j}v) \right|
& \leq & ||v_{x}||_{L^{\infty}(\mathbb{R}^{2})}
\left(\int(\partial^{j}w)^{2} \right)^{1/2}
\left(\int(\partial^{j}v)^{2} \right)^{1/2}\\
& \leq & c\left(\int[\,v_{x}^{2} + v_{xxx}^{2} +
v_{xy}^{2}\,] \right)^{1/2}||w||_{H^{|j|}(\mathbb{R}^{2})}\;
||v||_{H^{|j|}(\mathbb{R}^{2})}\\
& \leq & c\;||w||_{{Z}_{t}^{|j|}}\; ||v||_{{Z}_{t}^{|j|}}^{2}
\end{eqnarray*}
and
\begin{eqnarray*}
\left|\int w\,(\partial^{j}v_{x})\,(\partial^{j}v) \right|
& \leq &
\left|\int w_{x}\,(\partial^{j}v)^{2} \right|
\\
& \leq &
c\;||w_{x}||_{L^{\infty}(\mathbb{R}^{2})}\int
(\partial^{j}v)^{2} \\
& \leq & c\left(\int[\,w_{x}^{2} + w_{xxx}^{2} +
w_{xy}^{2}\,] \right)^{1/2}||v||_{{Z}_{t}^{|j|}}^{2}\\
& \leq & c\,||w||_{{Z}_{t}^{|j|}}\,||v||_{{Z}_{t}^{|j|}}^{2}.
\end{eqnarray*}
Therefore, we obtain
\begin{eqnarray*}
\partial_{t}\int(\partial^{j}v(\,\cdot\,,\,\cdot\,,\,\widetilde{t}))^{2} 
\leq c\,||w||_{{Z}_{t}^{|j|}}\,||v||_{{Z}_{t}^{|j|}}^{2}.
\end{eqnarray*}
Next, we take three $x$ derivatives of \eqref{e310}, multiply by
$2\,\partial^{j}v_{xxx}$ and integrate over $\mathbb{R}^{2}.$ Our
inequality becomes
\begin{align*}
\partial_{t}\int\partial^{j}v_{xxx}
& \leq  c\left|\int(\partial^{j}(w\,v_{x})_{xxx})
\,(\partial^{j}v_{xxx})  \right|\\
& \leq  c\left|\int
\partial^{j}(w_{xxx}\,v_{x} + 2\,w_{xx}\,v_{xx} + 2\,w_{x}\,v_{xxx} + w\,v_{xxxx})
(\partial^{j}v_{xxx})  \right|\\
& \leq  c\left|\int
\partial^{j}(w_{xxx}\,v_{x})\,(\partial^{j}v_{xxx}) \right|
+
c\left|\int\partial^{j}(w_{xx}\,v_{xx})\,
(\partial^{j}v_{xxx}) \right|\\
&  \qquad +\; c\left|\int\partial^{j}(w_{x}\,v_{xxx})\,
(\partial^{j}v_{xxx}) \right| + c
\left|\int\partial^{j}(w\,v_{xxxx})\,(\partial^{j}v_{xxx}) 
\right|\\
& \leq I_{1} + I_{2} + I_{3} + I_{4}.
\end{align*}
We will look at terms $I_{k},\;k=1,\,2,\,3,\,4$ below. For $I_{1}$
we have
\begin{eqnarray*}
\lefteqn{\left|\int
\partial^{j}(w_{xxx}\,v_{x})\,(\partial^{j}v_{xxx}) \right|}\\
& = & \left|\int
[\,\partial^{j}w_{xxx})\,v_{x}\,(\partial^{j}v_{xxx})\;+ \,\ldots
\,+\;
w_{xxx}\,(\partial^{j}v_{x})\,(\partial^{j}v_{xxx})\,] \right|\\
& \leq & ||v_{x}||_{L^{\infty}(\mathbb{R}^{2})}
\left(\int(\partial^{j}w_{xxx})^{2} \right)^{1/2}
\left(\int(\partial^{j}v_{xxx})^{2} \right)^{1/2}\\
&  & +\, \cdots\, +
||\partial^{j}v_{x}||_{L^{\infty}(\mathbb{R}^{2})}
\left(\int w_{xxx}^{2} \right)^{1/2}
\left(\int(\partial^{j}v_{xxx})^{2} \right)^{1/2}\\
& \leq & ||v||_{{Z}_{t}^{0}}\,||w||_{{Z}_{t}^{|j|}}\,||v||_{{Z}_{t}^{|j|}} +\,\cdots \, +
||v||_{{Z}_{t}^{|j|}}\,||w||_{{Z}_{t}^{0}}\,||v||_{{Z}_{t}^{|j|}}\\
& \leq & c\;||w||_{{Z}_{t}^{|j|}}\,||v||_{{Z}_{t}^{|j|}}^{2}.
\end{eqnarray*}
For $I_{2},$
\begin{eqnarray*}
\left|\int\partial^{j}(w_{xx}\,v_{xx})
\,(\partial^{j}v_{xxx}) \right| =
\left|\int[\,(\partial^{j}w_{xx})\,v_{xx}
+\,\ldots\,+
w_{xx}(\partial^{j}v_{xx})\,]\,(\partial^{j}v_{xxx}) \right|.
\end{eqnarray*}
To bound these terms, we will use the following anisotropic
imbedding in \cite{BIN}. For $2\leq n<6,$
\begin{eqnarray}
\label{e311}\left(\int_{\mathbb{R}^{2}}|u|^{n} \right)^{1/n}\leq
\left(\int_{\mathbb{R}^{2}}[\,u^{2} + u_{x}^{2} +
(\partial_{x}^{-1}u_{y})^{2}\,]\right)^{1/2}.
\end{eqnarray}
We will look at the most difficult terms to bound below.
\begin{eqnarray*}
\lefteqn{\left|\int
(\partial^{j}w_{xx})\,v_{xx}\,(\partial^{j}v_{xxx}) \right|}\\
& \leq &
\left(\int(\partial^{j}w_{xx})^{4} \right)^{1/4}
\left(\int(v_{xx})^{4} \right)^{1/4}
\left(\int(\partial^{j}v_{xxx})^{2} \right)^{1/2}\\
& \leq & \left(\int[\,(\partial^{j}\,w_{xx})^{2} +
(\partial^{j}\,w_{xxx})^{2} +
(\partial^{j}\,w_{xy})^{2}\,] \right)^{1/2}\\
&  & \times \left(\int[\,(v_{xx})^{2} +
(v_{xxx})^{2} + (v_{xy})^{2}\,] \right)^{1/2}
\left(\int(\partial^{j}v_{xxx})^{2} \right)^{1/2}\\
& \leq & ||w||_{{Z}_{t}^{|j|}}\,||v||_{{Z}_{t}^{|j|}}^{2},
\end{eqnarray*}
while,
\begin{eqnarray*}
\lefteqn{\left|\int w_{xx}\,(\partial^{j}v_{xx})
\,(\partial^{j}v_{xxx}) \right|}\\
& = &
c\left|\int w_{xxx}\,(\partial^{j}v_{xx})^{2} \right|\\
& \leq &
c\left(\int w_{xxx}^{2} \right)^{1/2}
\left(\int(\partial^{j}v_{xx})^{4} \right)^{1/2}\\
& \leq & c\;||w||_{{Z}_{t}^{0}}\left(\int
[\,(\partial^{j}v_{xx})^{2} + (\partial^{j}v_{xxx})^{2} +
(\partial^{j}v_{xy})^{2}\,] \right)^{1/2}\\
& \leq & c\;||w||_{{Z}_{t}^{|j|}}\,||v||_{{Z}_{t}^{|j|}}^{2}.
\end{eqnarray*}
For $I_{3},$
\begin{eqnarray*}
\lefteqn{\left|\int
\partial^{j}(w_{x}\,v_{xxx})\,(\partial^{j}v_{xxx}) \right|}\\
& = & \left|\int[\,(\partial^{j}w_{x})\,v_{xxx}
+\,\ldots\,+
w_{x}\,(\partial^{j}v_{xxx})\,]\,(\partial^{j}v_{xxx}) \right|\\
& \leq & ||\partial^{j}w_{x}||_{L^{\infty}(\mathbb{R}^{2})}
\left(\int v_{xxx}^{2} \right)^{1/2}
\left(\int(\partial^{j}v_{xxx})^{2} \right)^{1/2}\\
&  & +\,\ldots\,+ ||w_{x}||_{L^{\infty}(\mathbb{R}^{2})}
\left(\int(\partial^{j}v_{xxx})^{2} \right)^{1/2}\\
& \leq & c\,||w||_{{Z}_{t}^{|j|}}\,||v_{{Z}_{t}^{0}}\,||v||_{{Z}_{t}^{|j|}}
+\,\ldots\,+ c\,||w||_{{Z}_{t}^{0}}\,||v||_{{Z}_{t}^{|j|}}^{2}\\
& \leq & c\,||w||_{{Z}_{t}^{|j|}}\,||v||_{{Z}_{t}^{|j|}}^{2}.
\end{eqnarray*}
Lastly, for $I_{4},$
\begin{eqnarray*}
\left|\int\partial^{j}(w\,v_{xxxx})
\,(\partial^{j}v_{xxx}) \right| =
\left|\int[\,(\partial^{j}w)\,v_{xxxx} + \,\ldots\,
+
w\,(\partial^{j}v_{xxxx})\,]\,(\partial^{j}v_{xxx}) \right|.
\end{eqnarray*}
The first term is handled below. If $j=(0,\,0),$ then
\begin{eqnarray*}
\left|\int(\partial^{j}w)
\,v_{xxxx}\,(\partial^{j}v_{xxx}) \right| & = &
\left|\int w\,v_{xxxx}\,v_{xxx} \right|\\
& = & c\left|\int w_{x}\,v_{xxx}^{2} \right|\\
& \leq & c\,||w_{x}||_{L^{\infty}(\mathbb{R}^{2})}
\left(\int v_{xxx}^{2} \right)\\
& \leq & c\,||w||_{{Z}_{t}^{|j|}}\,||v||_{{Z}_{t}^{0}}^{2},
\end{eqnarray*}
while, if $|j|>0,$ then
\begin{eqnarray*}
\left|\int(\partial^{j}w)\,w_{xxxx}
\,(\partial^{j}v_{xxx}) \right| & \leq &
||\partial^{j}w||_{L^{\infty}(\mathbb{R}^{2})}
\left(\int v_{xxxx}^{2}\right)^{1/2}
\left(\int(\partial^{j}v_{xxx})^{2}\right)^{1/2}\\
& \leq & ||w||_{{Z}_{t}^{|j|}}\,||v||_{{Z}_{t}^{|j|}}^{2}.
\end{eqnarray*}
The last term in $I_{4}$ is handled below
\begin{eqnarray*}
\left|\int w\,(\partial^{j}v_{xxxx})
\,(\partial^{j}v_{xxx}) \right| & = &
c\left|\int w_{x}\,(\partial^{j}v_{xxx})^{2} \right|\\
& \leq & c\,||w_{x}||_{L^{\infty}(\mathbb{R}^{2})}
\left(\int(\partial^{j}v_{xxx})^{2} \right)\\
& \leq & c\,||w||_{{Z}_{t}^{0}}\,||v||_{{Z}_{t}^{|j|}}^{2}\\
& \leq & c\,||w||_{{Z}_{t}^{|j|}}\,||v||_{{Z}_{t}^{|j|}}^{2}.
\end{eqnarray*}
Consequently, we conclude
\begin{eqnarray*}
\partial_{t}\int(\partial^{j}v_{xxx})^{2} \leq
c\,||w||_{{Z}_{t}^{|j|}}\,||v||_{{Z}_{t}^{|j|}}^{2}.
\end{eqnarray*}
Next we take two $y$ derivatives of \eqref{e310}, multiply by
$2\,(\partial^{j}v_{yy}),$ and integrate over $\mathbb{R}^{2}.$
Therefore, we have
\begin{align*}
\partial_{t}\int(\partial^{j}v_{yy})^{2} 
& \leq c\,\left|\int\partial^{j}(w\,v_{x})_{yy}
\,(\partial^{j}v_{yy}) \right|\\
& \leq c\,\left|\int\partial^{j}(w_{yy}\,v_{x} +
2\,w_{y}\,v_{xy} + w\,v_{xyy})
\,(\partial^{j}v_{yy}) \right|\\
& \leq c\,\left|\int\partial^{j}(w_{yy}\,v_{x})\,(\partial^{j}v_{yy}) \right|
+ c\left|\int\partial^{j}(w_{y}\,v_{xy})\,(\partial^{j}v_{yy}) \right| \\
& \qquad + c\left|\int\partial^{j}(w\,v_{xyy}) \,(\partial^{j}v_{yy}) \right|\\
& \leq I_{5} + I_{6} + I_{7}.
\end{align*}
First, we look at $I_{5},$
\begin{eqnarray*}
\lefteqn{\left|\int\partial^{j}(w_{yy}\,v_{x})
\,(\partial^{j}v_{yy}) \right|}\\
& = & \left|\int[\,(\partial^{j}w_{yy})\,v_{x} +
\,\ldots\, +
w_{yy}\,(\partial^{j}v_{x})\,]\,(\partial^{j}v_{yy}) \right|\\
& \leq & c\,||v_{x}||_{L^{\infty}(\mathbb{R}^{2})}
\left(\int(\partial^{j}w_{yy})^{2} \right)^{1/2}
\left(\int(\partial^{j}v_{yy})^{2} \right)^{1/2}\\
&  & +\,\ldots\,+
c\,||\partial^{j}v_{x}||_{L^{\infty}(\mathbb{R}^{2})}
\left(\int w_{yy}^{2} \right)^{1/2}
\left(\int(\partial^{j}v_{yy})^{2} \right)^{1/2}\\
& \leq & c\,||v||_{{Z}_{t}^{0}}\,||w||_{{Z}_{t}^{|j|}}\,||v||_{Z_{t}^{|j|}}
+\,\ldots\,+ c\,||v||_{{Z}_{t}^{|j|}}\,||w||_{{Z}_{t}^{0}}\\
& \leq & c\,||w||_{{Z}_{t}^{|j|}}\,||v||_{{Z}_{t}^{|j|}}^{2}.
\end{eqnarray*}
For $I_{6},$
\begin{eqnarray*}
\lefteqn{\left|\int\partial^{j}(w_{y}\,v_{xy})
\,(\partial^{j}v_{yy}) \right|}\\
& = & \left|\int[\,(\partial^{j}w_{y})\,v_{xy} +
\,\ldots\, +
w_{y}\,(\partial^{j}v_{xy})\,]\,(\partial^{j}v_{yy}) \right|\\
& \leq & c\,||w_{y}||_{L^{\infty}(\mathbb{R}^{2})}
\left(\int(\partial^{j}v_{xy})^{2} \right)^{1/2}
\left(\int(\partial^{j}v_{yy})^{2} \right)^{1/2}\\
&  & +\,\ldots\,+
c\,||\partial^{j}w_{y}||_{L^{\infty}(\mathbb{R}^{2})}
\left(\int v_{xy}^{2} \right)^{1/2}
\left(\int(\partial^{j}v_{yy})^{2} \right)^{1/2}\\
& \leq & c\,||w||_{{Z}_{t}^{|j|}}\,||v||_{{Z}_{t}^{0}}\,||v||_{{Z}_{t}^{|j|}}
+\,\ldots\,+ c\,||w||_{{Z}_{t}^{0}}\,||v||_{{Z}_{t}^{0}}\\
& \leq & c\,||w||_{{Z}_{t}^{|j|}}\,||v||_{{Z}_{t}^{|j|}}^{2}.
\end{eqnarray*}
Lastly, for $I_{7},$
\begin{eqnarray*}
\left|\int\partial^{j}(w\,v_{xyy})\,(\partial^{j}v_{yy}) \right|
= \left|\int[\,(\partial^{j}w)\,v_{xyy} +
\,\ldots\,+
w\,(\partial^{j}v_{xyy})\,]\,(\partial^{j}v_{yy}) \right|.
\end{eqnarray*}
We will look at the first and last of these terms below. The rest of
these terms are handled similarly. For the first term, if
$j=(0,\,0),$ then we have
\begin{eqnarray*}
\left|\int(\partial^{j}w)
\,v_{xyy}\,(\partial^{j}v_{yy}) \right| & = &
\left|\int w\,v_{xyy}\,v_{yy} \right|\\
& = & c\left|\int w_{x}\,v_{yy}^{2} \right|\\
& \leq & c\,||w_{x}||_{L^{\infty}(\mathbb{R}^{2})}
\int v_{yy}^{2} \\
& \leq & c\,||w||_{{Z}_{t}^{0}}\,||v||_{{Z}_{t}^{|j|}}^{2},
\end{eqnarray*}
while for $|j|>0,$
\begin{eqnarray*}
\left|\int
(\partial^{j}w)\,v_{xyy}\,(\partial^{j}v_{yy}) \right| & \leq
& ||\partial^{j}w||_{L^{\infty}(\mathbb{R}^{2})}
\left(\int v_{xyy}^{2}\right)^{1/2}
\left(\int(\partial^{j}v_{yy})^{2}\right)^{1/2}\\
& \leq & ||w||_{{Z}_{t}^{|j|}}\,||v||_{{Z}_{t}^{|j|}}^{2}.
\end{eqnarray*}
The last term for $I_{7}$ is bounded as follows,
\begin{eqnarray*}
\left|\int w\,(\partial^{j}v_{xyy})\,(\partial^{j}v_{yy}) \right|
& = &
c\left|\int w_{x}\,(\partial^{j}v_{yy})^{2} \right|\\
& \leq &
c\,||w_{x}||_{L^{\infty}(\mathbb{R}^{2})}\int(\partial^{j}v_{yy})^{2} \\
& \leq & c\,||w||_{{Z}_{t}^{|j|}}\,||v||_{{Z}_{t}^{|j|}}^{2}.
\end{eqnarray*}
Now apply one $t$ derivative to \eqref{e310}, multiply by
$2\,(\partial^{j}v_{t})$ and integrate over $\mathbb{R}^{2}.$ We
arrive at the following inequality,
\begin{eqnarray*}
\partial_{t}\int(\partial^{j}v_{t})  & \leq
& c\left|\int(\partial^{j}(w\,v_{x})_{t})
\,(\partial^{j}v_{t}) \right|\\
& \leq &
c\left|\int\partial^{j}(w_{t}\,v_{x})\,(\partial^{j}v_{t}) \right|
+
c\left|\int\partial^{j}(w\,v_{xt})\,(\partial^{j}v_{t})\,dx\,dt\right|\\
& = & I_{9} + I_{10}.
\end{eqnarray*}
For $I_{9},$ we have
\begin{eqnarray*}
\lefteqn{\left|\int\partial^{j}(w_{t}\,v_{x})
\,(\partial^{j}v_{t}) \right|}\\
& \leq &
c\left|\int(\partial^{j}w_{t})\,v_{x}\,(\partial^{j}v_{t}) \right|
+\,\ldots\,+
c\left|\int w_{t}\,(\partial^{j}v_{x})\,(\partial^{j}v_{t}) \right|\\
& \leq & c\,||v_{x}||_{L^{\infty}(\mathbb{R}^{2})}
\left(\int(\partial^{j}w_{t})^{2} \right)^{1/2}
\left(\int(\partial^{j}v_{t})^{2} \right)^{1/2}\\
&  & +\,\ldots\,+ ||\partial^{j}v_{x}||_{L^{\infty}(\mathbb{R}^{2})}
\left(\int(w_{t})^{2} \right)^{1/2}
\left(\int(\partial^{j}v_{t})^{2} \right)^{1/2}\\
& \leq & c\,||v||_{{Z}_{t}^{0}}\,||w||_{Z_{t}^{|j|}}\,||v||_{{Z}_{t}^{|j|}}
+\,\ldots\,+
c\,||v||_{{Z}_{t}^{|j|}}\,||w||_{{Z}_{t}^{0}}\,||v||_{{Z}_{t}^{|j|}}\\
& \leq & c\,||w||_{{Z}_{t}^{|j|}}\,||v||_{{Z}_{t}^{|j|}}^{2}.
\end{eqnarray*}
Next we look at $I_{10}.$ If $j=(0,\,0),$ we have
\begin{eqnarray*}
\left|\int w\,v_{xt}\,v_{t} \right| & = &
\left|\int w_{x}\,v_{t}^{2} \right|\\
& \leq &
||w_{x}||_{L^{\infty}(\mathbb{R}^{2})}\int v_{t}^{2} \\
& \leq & c\,||w||_{{Z}_{t}^{|j|}}\,||v||_{{Z}_{t}^{|j|}}^{2}.
\end{eqnarray*}
If $j\neq (0,\,0),$ we have
\begin{align*}
\left|\int\partial^{j}(w\,v_{xt})\,(\partial^{j}v_{t}) \right|
& =
c\left|\int(\partial^{j}w)\,v_{xt}\,(\partial^{j}v_{t}) \right|
+\,\ldots \\
& \qquad +
c\left|\int w\,(\partial^{j}v_{xt})\,(\partial^{j}v_{t}) \right|\\
& = I_{10}(a) +\,\ldots\,+ I_{10}(\widetilde{a}).
\end{align*}
Now for $I_{10}(a),$ we use the following estimate
\begin{eqnarray*}
\left|\int(\partial^{j}w)\,v_{xt}\,(\partial^{j}v_{t}) \right|
& \leq & c\,||\partial^{j}w||_{L^{\infty}(\mathbb{R}^{2})}
\left(\int v_{xt}^{2} \right)^{1/2}
\left(\int(\partial^{j}v_{t})^{2} \right)^{1/2}\\
& \leq & c\,||w||_{{Z}_{t}^{|j|}}\,||v||_{{Z}_{t}^{|j|}}^{2}.
\end{eqnarray*}
While for $I_{10}(\widetilde{a}),$ we use the following estimate
\begin{eqnarray*}
\left|\int w\,(\partial^{j}v_{xt})\,(\partial^{j}v_{t}) \right|
& = &
c\left|\int w_{x}\,(\partial^{j}v_{t})^{2} \right|\\
& \leq &
c\,||w_{x}||_{L^{\infty}(\mathbb{R}^{2})}\int(\partial^{j}v_{t})^{2} \\
& \leq & c\,||w||_{{Z}_{t}^{|j|}}\,||v||_{{Z}_{t}^{|j|}}^{2}.
\end{eqnarray*}
Therefore, for $0\leq\widetilde{t}\leq t,$ we conclude that
\begin{align*}
& \partial_{t}\int[\,(\partial^{j}v(\,\cdot\,,\,\cdot\,,\,\widetilde{t}))^{2}
+ (\partial^{j}v_{xxx}(\,\cdot\,,\,\cdot\,,\,\widetilde{t}))^{2} +
(\partial^{j}v_{yy}(\,\cdot\,,\,\cdot\,,\,\widetilde{t}))^{2} +
(\partial^{j}v_{t}(\,\cdot\,,\,\cdot\,,\,\widetilde{t}))^{2}\,] \\
& \qquad \leq
c\,||w||_{{Z}_{t}^{|j|}}\,||v||_{{Z}_{t}^{|j|}}^{2}.
\end{align*}
Integrating with respect to $t,$ we obtain
\begin{eqnarray*}
||v||_{Z_{t}^{|j|}}^{2}\leq ||v(\,\cdot\,,\,\cdot\,,\,0))||_{H^{(|j|
+ 3,\,|j| + 2)}(\mathbb{R}^{2})}^{2} +
||v(\,\cdot\,,\,\cdot\,,\,0))||_{H^{|j|}(\mathbb{R}^{2})}^{2} +
c\,t\,||w||_{Z_{t}^{|j|}}\,||v||_{Z_{t}^{|j|}}^{2},
\end{eqnarray*}
as desired. \hfill $\square$

\renewcommand{\theequation}{\thesection.\arabic{equation}}
\setcounter{equation}{0}\section{Uniqueness and Existence of a
local solution} In this section, we will prove that for $\phi\in
X^{N}(\mathbb{R}^{2})$ there exists a unique solution of
\eqref{e201}-\eqref{e202} in
$L^{\infty}([0,\,T]:\,X^{N}(\mathbb{R}^{2})),$ where the time $T$
depends only $||\phi||_{X^{0}(\mathbb{R}^{2})}.$  First we prove uniqueness of
solutions.\\
\\
{\bf Theorem 4.1} (Uniqueness). {\it Let $\phi\in
X^{0}(\mathbb{R}^{2})$ and $0<T<+\infty.$ Then there is at most
one solution of \eqref{e201}-\eqref{e202} in
$L^{\infty}([0,\,T]:\,X^{0}(\mathbb{R}^{2}))$ with initial data
$u(x,\,y,\,0)=\phi(x,\,y).$}\\
\\
{\it Proof.} Assume that $u,\,$ $v\in
L^{\infty}([0,\,T]:\,X^{0}(\mathbb{R}^{2}))$ are two solutions of
\eqref{e201}-\eqref{e202} with $u_{t},\,$ $v_{t}\in
L^{\infty}([0,\,T]:\,L^{2}(\mathbb{R}^{2})),$ so all integrations
below are justified and with the same initial data, in fact, with
$(u - v)(x,\,y,\,0)=0.$ Then
\begin{eqnarray}
\label{e401}&  & (u - v)_{t} + (u - v)_{xxx} + (u - v)_{x} +
(u\,u_{x} - v\,v_{x}) -
\partial_{x}^{-1}(u - v)_{yy} =0.
\end{eqnarray}
By \eqref{e401},
\begin{eqnarray}
\label{e402}&  & (u - v)_{t} + (u - v)_{xxx} + (u - v)_{x} + (u -
v)\,u_{x} + (u - v)_{x}\,v -
\partial_{x}^{-1}(u - v)_{yy} =0.
\end{eqnarray}
Multiplying \eqref{e402} by $2\,(u - v)$ and integrating with
respect to $(x,\,y)$ over $\mathbb{R}^{2},$
\begin{eqnarray}
\label{e403}&  & 2\int (u - v)\,(u -
v)_{t} + 2\int (u - v)\,(u -
v)_{xxx} +
2\int (u - v)\,(u - v)_{x}   \\
&  & +\;2\int (u - v)^{2}\,u_{x} +
2\int (u - v)\,(u - v)_{x}\,v  -
2\int (u - v)\,\partial_{x}^{-1}(u -
v)_{yy} =0.\nonumber
\end{eqnarray}
Integrating by parts each term in \eqref{e403} we obtain
\begin{eqnarray}
\partial_{t}\int (u - v)^{2} & = &
-\;2\int (u - v)^{2}\,u_{x} +
\int (u - v)^{2}\,v_{x} \nonumber \\
& \leq & c\left(||u_{x}||_{L^{\infty}(\mathbb{R}^{2})} +
||v_{x}||_{L^{\infty}(\mathbb{R}^{2})}\right)\int(u
- v)^{2} \nonumber \\
\label{e404}& \leq & c\left(\,||u||_{X^{0}(\mathbb{R}^{2})} +
||v||_{X^{0}(\mathbb{R}^{2})}\,\right)\int (u - v)^{2}
\end{eqnarray}
Using Gronwall's inequality and the fact that $(u - v)$ vanishes
at $t=0,$ it follows that $u=v.$ This proves the uniqueness of the
solution. \hfill $\square$ \\
\\
Now we consider existence of solutions to \eqref{e201}-\eqref{e202}.  Our plan is to show
that for $\phi \in X^N(\mathbb R^2)$ there exists a solution $u \in L^\infty([0,T]:X^N(\mathbb R^2))$
for a time $T$ depending only on $||\phi||_{X^0}$.  In order to prove this we must first prove
a preliminary result by introducing the following function space.  Let
\begin{eqnarray}
\label{e301}Y^{N}(\mathbb{R}^{2}) =
\left\{u:\;u,\,u_{xxx},\,u_{yy},\,
\frac{\eta^{2}}{\xi}\,\widehat{u}\in H^{N}(\mathbb{R}^{2})\right\}
\end{eqnarray}
with the accompanying norm
\begin{eqnarray}
\label{e302}||u||_{Y^{N}(\mathbb{R}^{2})}^{2}=\int_{\mathbb{R}^{2}}\left(u^{2}
+ \sum_{|j|\leq N}[\,(\partial^{j}u_{xxx})^{2} +
(\partial_{x}^{-1}\partial^{j}u_{yy})^{2} +
(\partial^{j}u_{yy})^{2}\,]\right)dx\,dy
\end{eqnarray}
where $j=(\alpha_1,\,\alpha_2)$ and $|j| = \alpha_1 + \alpha_2.$
We will begin by showing that, for $\phi \in Y^N(\mathbb R^2)$, there exists a solution
$u$ of \eqref{e201} such that $u \in L^\infty([0,T]: Y^N(\mathbb R^2))$ for a time $T$
depending only on $||\phi||_{Y^0}$.  Then we will prove a differential inequality of
the form
\[
\partial_t \left(\int u^2 + u_{xxx}^2 + (\partial_x^{-1}u_{yy})^2\right) \leq
\left(\int u^2 + u_{xxx}^2 + (\partial_x^{-1}u_{yy})^2\right)^{3/2},
\]
to show that in fact the solution $u$ obtained in Theorem 4.2 is in $L^\infty([0,T'];
X^0(\mathbb R^2))$ for a time $T'$ depending only on $||\phi||_{X^0(\mathbb R^2)}$.
With these ideas in mind we state our existence theorem.
\\

\nd{\bf Theorem 4.2} (Existence). {\it Let $k_{0}>0$ and $N$ be an
integer $\geq 0.$ Then there exists a time $0<T<+\infty$ depending
only on $k_{0}$ such that for all $\phi\in Y^{N}(\mathbb{R}^{2})$
with $||\phi||_{Y^{0}(\mathbb{R}^{2})}\leq k_{0}$ there exists a
solution of \eqref{e201}, $u\in
L^{\infty}\left([0,\,T]:\,Y^{N}(\mathbb{R}^{2})\right)$ such that
$u(x,\,y,\,0)=\phi(x,\,y).$}\\
\\

The method of proof is as follows. As discussed in section 3, we begin by approximating
\eqref{e201} by the linear equation \eqref{e305}. We construct the mapping
\begin{eqnarray*}
\Pi : Z_{T}^{N} \rightarrow Z_{T}^{N}
\end{eqnarray*}
where the initial condition is given by $u^{(n)}(x,\,y,\,0)=\phi(x,\,y)$
and the first approximation is given by $u^{(0)}(x,y,t) = \phi(x,y)$.  Subsequent
approximations are given by $u^{(n)} = \Pi(u^{(n-1)})$ for $n \geq 1$.
Equation \eqref{e305} is a linear equation which by Lemma 3.1 can
be solved at each iteration. We show that the sequence of solutions $\{u^{(n)}\}$
to our linear equation is bounded in $L^\infty([0,T];Y^0(\mathbb R^2))$ for a time
$T$ depending only on $||\phi||_{Y^0}$.  We then show that there is a subsequence of
solutions to our approximate equations which converges to a solution $u \in L^\infty([0,T];
Y^0(\mathbb R^2))$ of \eqref{e201}.  Lastly, we show that if $\phi \in Y^N(\mathbb R^2)$
for $N > 0$, then our solution $u \in L^\infty([0,T];Y^N(\mathbb R^2))$ where the time
$T$ depends only on $||\phi||_{Y^0}$.\\
\\

\nd{\it Proof.} It suffices to prove this result for
$\phi\in\bigcap_{N\geq 0}H^{N}(\mathbb{R}^{2})$ and
$\partial_{x}^{-1}\phi_{yy}\in\bigcap_{N\geq
0}H^{N}(\mathbb{R}^{2}).$ We can then use the same approximation
procedure as before to prove the result for general initial data.
Let $u^{(n)}$ be a solution of \eqref{e301} with initial data $u^{(n)}(x,y,0) =
\phi(x,y)$ and where the first approximation is given by $u^{(0)}(x,y,t) = \phi(x,y)$.
By Lemma 3.2, we know that
\begin{eqnarray}
\label{e405}||u^{(n)}||_{{Z}_{t}^{0}}^{2} \leq
||u^{(n)}(\,\cdot\,,\,\cdot\,,\,0)||_{H^{(3,\,2)}(\mathbb{R}^{2})}^{2}
+
||u_{t}^{(n)}(\,\cdot\,,\,\cdot\,,\,0)||_{L^{2}(\mathbb{R}^{2})}^{2}
+ c\,t\,||u^{(n - 1)}||_{{Z}_{t}^{0}}\;||u^{(n)}||_{{Z}_{t}^{0}}^{2}.
\end{eqnarray}
Further, using the fact that $||\phi||_{Y^{0}}\leq k_{0}$, we have
\begin{eqnarray*}
\lefteqn{||u^{(n)}(\,\cdot\,,\,\cdot\,,\,0)||_{H^{(3,\,2)}(\mathbb{R}^{2})}^{2}
+
||u_{t}^{(n)}(\,\cdot\,,\,\cdot\,,\,0)||_{L^{2}(\mathbb{R}^{2})}^{2}}
\\
& = &
||u^{(n)}(\,\cdot\,,\,\cdot\,,\,0)||_{H^{(3,\,2)}(\mathbb{R}^{2})}^{2}\\
&  & +
\int [\,u_{xxx}^{(n)}(\,\cdot\,,\,\cdot\,,\,0) +
u_{x}^{(n)}(\,\cdot\,,\,\cdot\,,\,0) -
\partial_{x}^{-1}u_{yy}^{(n)}(\,\cdot\,,\,\cdot\,,\,0) +
u^{(n -
1)}(\,\cdot\,,\,\cdot\,,\,0)u_{x}^{(n)}(\,\cdot\,,\,\cdot\,,\,0)\,]^{2} \\
& \leq & ||\phi||_{Y^{0}(\mathbb{R}^{2})}^{2} +
c\int [\,\phi_{xxx}^{2} + \phi_{x}^{2} +
(\partial_{x}^{-1}\phi_{yy})^{2} + (\phi\,\phi_{x})^{2}\,]
\\
& \leq & C\,||\phi||_{Y^{0}(\mathbb{R}^{2})}^{2} \leq C\,k_{0}^{2},
\end{eqnarray*}
where $\;C\,$ is independent of $n.\;$ Define
$c_{0}=\left(\frac{C}{2}k_{0}^{2} + 1\right).\;$ Let $T_{0}^{(n)}$
be the maximum time such that $||u^{(j)}||_{Z_{t}^{0}}\leq c_{0}$
for $0\leq t\leq T_{0}^{(n)},$ $0\leq j\leq n.$ That is
\begin{eqnarray*}
T_{0}^{(n)}=\sup\{t\in [0,\,T_{0}^{(n)}]:\;||u^{(j)}||_{{Z}_{t}^{0}}\leq c_{0}\quad\mbox{for}\quad 0\leq j\leq n\}.
\end{eqnarray*}
Therefore,
\begin{eqnarray}
\label{e406}||u^{(n)}||_{{Z}_{t}^{0}}^{2} & \leq &
||u^{(n)}(\,\cdot\,,\,\cdot\,,\,0)||_{H^{(3,\,2)}(\mathbb{R}^{2})}^{2}
+
||u_{t}^{(n)}(\,\cdot\,,\,\cdot\,,\,0)||_{L^{2}(\mathbb{R}^{2})}^{2}
+ c\,t\,||u^{(n - 1)}||_{{Z}_{t}^{0}}\;||u^{(n)}||_{{Z}_{t}^{0}}^{2} \nonumber \\
& \leq & C\,k_{0}^{2} + c\,t\,c_{0}^{3}.
\end{eqnarray}
{\it Claim:} $T_{0}^{(n)}$ does not approach $0.$\\
\\
On the contrary, assume that $T_{0}^{(n)}\rightarrow 0.$ Since
$||u^{(n)}(\,\cdot\,,\,\cdot\,,\,t)||_{{Z}_{t}^{0}}$ is
continuous for $t\geq 0,$ there exists $\tau\in [0,\,T]$ such that
$c_{0}=||u^{(j)}(\,\cdot\,,\,\cdot\,,\,\tau)||_{{Z}_{\tau}^{0}}$ for $0\leq\tau\leq T_{0}^{(n)},$ $0\leq j\leq n.$
Then, by \eqref{e406} we have
\begin{eqnarray}
\label{e407}c_{0}^{2} \leq C\,k_{0}^{2} +
c\,T_{0}^{(n)}\,c_{0}^{3}.
\end{eqnarray}
As $n\rightarrow \infty,$ we have
\begin{eqnarray}
\label{e408}\left(\frac{C}{2}\,k_{0}^{2} + 1\right)^{2} \leq
C\,k_{0}^{2}\quad \Longrightarrow \quad \frac{C^{2}}{4}\,k_{0}^{4}
+ 1 \leq 0
\end{eqnarray}
which is a contradiction. Consequently $T_{0}^{(n)}\not\rightarrow
0.$ Choosing $T=T(c_0)$ sufficiently small, and $T$ not depending on
$n,$ one concludes that
\begin{eqnarray}
\label{e409}||u^{(n)}||_{{Z}_{t}^{0}}^{2} \leq
c\quad\mbox{for}\quad 0\leq t\leq T.
\end{eqnarray}
This show that $T_{0}^{(n)}\geq T.$ Hence from \eqref{e409} we see
that there exists a bounded sequence of solutions $u^{(n)}\in {Z}_{T}^{0}$ and therefore a subsequence $u^{(n_{j})}\equiv u^{(n)}$
such that
\begin{eqnarray*}
u^{(n)}\stackrel{*}\rightharpoonup u\quad \mbox{weakly in}\quad
L^{\infty}([0,\,T]:\,H^{(3,\,2)}(\mathbb{R}^{2}))
\end{eqnarray*}
\begin{eqnarray*}
u_{t}^{(n)}\stackrel{*}\rightharpoonup u_{t}\quad \mbox{weakly
in}\quad L^{\infty}([0,\,T]:\,L^{2}(\mathbb{R}^{2})).
\end{eqnarray*}
Therefore, by Lions-Aubin's compactness theorem there is a
subsequence $u^{(n_{j})}\equiv u^{(n)}$ such that
$u^{(n)}\rightarrow u$ strongly on
$L^{\infty}([0,\,T]:\,H_{loc}^{1}(\mathbb{R}^{2})).$ Now it
remains to show that each term in \eqref{e305} converges to its
correct limit. First, $u_{xxx}^{(n)}\stackrel{*}\rightharpoonup
u_{xxx}$ weakly on $L^{\infty}([0,\,T]:\,L^{2}(\mathbb{R}^{2})).$
Similarly $u_{t}^{(n)}\rightarrow u_{t}$ and
$u_{x}^{(n)}\rightarrow u_{x}$ weak$^{*}$ in
$L^{\infty}([0,\,T]:\,L^{2}(\mathbb{R}^{2})).$ Now we will show
that the nonlinear term converges to its correct limit. First,
$u^{(n - 1)}\rightarrow u$ strongly in
$L^{\infty}([0,\,T]:\,H_{loc}^{1}(\mathbb{R}^{2})).$ Moreover,
$u_{x}^{(n - 1)}\stackrel{*}\rightharpoonup $ weakly in
$L^{\infty}([0,\,T]:\,L^{2}(\mathbb{R}^{2})).$ Therefore,
\begin{eqnarray*}
u^{(n - 1)}\,u_{x}^{(n)}\stackrel{*}\rightharpoonup u\,u_{x}\quad
\mbox{weakly in}\quad L^{\infty}([0,\,T]:\,L^{2}(\mathbb{R}^{2})).
\end{eqnarray*}
Consequently,
\begin{equation}
\begin{split}
\partial_{x}^{-1}u_{yy}^{(n)}
& = u_{t}^{(n)} + u_{xxx}^{(n)} + u_{x}^{(n)} +
u^{(n)}\,u_{x}^{(n)}\\
&  \stackrel{*}\rightharpoonup u_{t} + u_{xxx} +
u_{x} + u\,u_{x}\quad \mbox{weakly in} \quad
L^{\infty}([0,\,T]:\,L^{2}(\mathbb{R}^{2})).
\end{split}
\end{equation}
But, also note that
\begin{eqnarray*}
u_{yy}^{(n)}\stackrel{*}\rightharpoonup u_{yy}\quad \mbox{weakly
in}\quad L^{\infty}([0,\,T]:\,L^{2}(\mathbb{R}^{2})).
\end{eqnarray*}
Therefore
\begin{eqnarray*}
\partial_{x}^{-1}u_{yy}^{(n)}
\stackrel{*}\rightharpoonup \partial_{x}^{-1}u_{yy}\quad
\mbox{weakly in}\quad L^{\infty}([0,\,T]:\,L^{2}(\mathbb{R}^{2}))
\end{eqnarray*}
and consequently $u$ is a solution to \eqref{e201}. Now, we prove
that there exists a solution to \eqref{e201} with $u\in
L^{\infty}([0,\,T]:\,Y^{N}(\mathbb{R}^{2}))$ for the time $T$ chosen
above. We already know that there is a solution $u\in
L^{\infty}([0,\,T]:\,Y^{0}(\mathbb{R}^{2})).$ Therefore, it
suffices to show that the approximating sequence $u^{(n)}$ is
bounded in ${Z}_{T}^{N}$ and thus, by the convergence arguments
above, our solution $u$ is in
$L^{\infty}([0,\,T]:\;Y^{N}(\mathbb{R}^{2})).$ Again, by Lemma 3.1,
we know our linearized equation can be solved in any interval of
time in which the coefficients are defined. Therefore, for each
iterate, $||u^{(n)}||_{{Z}_{t}^{N}}$ is continuous in $t\in
[0,\,T].$ By Lemma 3.2, it follows that
\begin{eqnarray}
||u^{(n)}||_{{Z}_{t}^{N}}^{2}  & \leq &
||u^{(n)}(\,\cdot\,,\,\cdot\,,\,0)||_{H^{(N + 3,\,N +
2)}(\mathbb{R}^{2})}^{2} \nonumber \\
\label{e410}&  & +\;
||u_{t}^{(n)}(\,\cdot\,,\,\cdot\,,\,0)||_{H^{N}(\mathbb{R}^{2})}^{2}
+ c\,t\,||u^{(n - 1)}||_{{Z}_{t}^{N}}\;||u^{(n)}||_{{Z}_{t}^{N}}^{2}.
\end{eqnarray}
On the other hand, as before and using $||\phi||_{Y^{N}}\leq
k_{N}$ we obtain
\begin{eqnarray*}
||u^{(n)}(\,\cdot\,,\,\cdot\,,\,0)||_{H^{(N + 3,\,N +
2)}(\mathbb{R}^{2})}^{2} +
||u_{t}^{(n)}(\,\cdot\,,\,\cdot\,,\,0)||_{H^{N}(\mathbb{R}^{2})}^{2}
&\leq  C\,k_{N}^{2},
\end{eqnarray*}
where $k_{N}$ is independent of $n.$ Define
$c_{N}=\left(\frac{C}{2}\,k_{N}^{2} + 1\right).$ Let $T_{N}^{(n)}$
be the largest time that $||u^{(j)}||_{{Z}_{t}^{N}}\leq c_{N}$
for $0\leq t\leq T_{N}^{(n)},$ $0\leq j\leq n.$ That is,
\begin{eqnarray*}
T_{N}^{(n)}=\sup\{t\in [0,\,T_{N}^{(n)}]:\;||u^{(j)}||_{{Z}_{t}^{N}}\leq c_{N}\quad\mbox{for}\quad 0\leq j\leq n\}.
\end{eqnarray*}
Therefore, for $0\leq t\leq T_{N}^{(n)},$
\begin{eqnarray}
\label{e411}||u^{(n)}||_{{Z}_{t}^{N}}^{2}\leq C\,k_{N}^{2} +
c\,t\,c_{N}^{3}.
\end{eqnarray}
{\it Claim:} $T_{N}^{(n)}$ does not approach $0.$\\
\\
On the contrary, assume that $T_{N}^{(n)}\rightarrow 0.$ Since
$||u^{(n)}(\,\cdot\,,\,\cdot\,,\,t)||_{{Z}_{t}^{N}}$ is
continuous for $t\geq 0,$ there exists $\tau\in [0,\,T]$ such that
$c_{N}=||u^{(j)}(\,\cdot\,,\,\cdot\,,\,\tau)||_{{Z}_{\tau}^{N}}$ for $0\leq\tau\leq T_{N}^{(n)},$ $0\leq j\leq n.$
Then, by \eqref{e411} we have
\begin{eqnarray}
\label{e412}c_{N}^{2} \leq C\,k_{N}^{2} +
c\,T_{N}^{(n)}\,c_{N}^{3}.
\end{eqnarray}
As $n\rightarrow \infty,$ we have
\begin{eqnarray}
\label{e413}\left(\frac{C}{2}\,k_{N}^{2} + 1\right)^{2} \leq
C\,k_{N}^{2}\quad \Longrightarrow \quad \frac{C^{2}}{4}\,k_{N}^{4}
+ 1 \leq 0
\end{eqnarray}
which is a contradiction. Consequently $T_{N}^{(n)}\not\rightarrow
0.$ Choosing $T_N$ sufficiently small, and $T_N$ not depending on
$n,$ one concludes that
\begin{eqnarray}
\label{e414}||u^{(n)}||_{{Z}_{t}^{N}}^{2} \leq
c\quad\mbox{for}\quad 0\leq t\leq T_N.
\end{eqnarray}
This show that $T_{N}^{(n)}\geq T_N.$ Now, let
\begin{eqnarray*}
T_{N}^{*} =\sup\{t\in [0,\,T_{N}^{*}]:\;u\in {Z}_{t}^{N}\}.
\end{eqnarray*}
We claim that $T_{N}^{*}\geq T$ and therefore, a time of existence
can be chosen depending only on $||\phi||_{Y^{0}}.$ By Lemma 3.1
the linear equation \eqref{e305} can be solved in any interval of
time in which the coefficients are defined, and thus
$T_{N}^{*}\geq T.$ \hfill $\square$ \\
\\
Now we want to improve our existence theorem. In particular, we
want to show that the solution $u\in
L^{\infty}([0,\,T]:\,Y^{N}(\mathbb{R}^{2}))$ found in Theorem 4.2
is in $L^{\infty}([0,\,T']:\,X^{N}(\mathbb{R}^{2}))$ for a time
$T'$ depending only on $||\phi||_{X^{0}(\mathbb{R}^{2})}.$  In order to do so,
we first prove a differential inequality. \\
\\
{\bf Lemma 4.3.} {\it Let $u$ be the solution to our main
equation in $L^{\infty}([0,\,T]:\,Y^{N}(\mathbb{R}^{2})).$ Then
for any $0\leq t\leq T,$ we have}
\begin{eqnarray}
\lefteqn{\partial_{t}\int_{\mathbb{R}^{2}}\left(u^{2} +
\sum_{|j|\leq N}[\,(\partial^{j}u_{xxx})^{2} +
(\partial^{j}(\partial_{x}^{-1}u_{yy}))^{2}\,]\right)dx\,dy}\nonumber
\\
\label{e415}& \leq & c\left[\int_{\mathbb{R}^{2}}\left(u^{2} +
\sum_{|j|\leq N}[\,(\partial^{j}u_{xxx})^{2} +
(\partial^{j}(\partial_{x}^{-1}u_{yy}))^{2}\,]\right)dx\,dy\right]^{3/2}.
\end{eqnarray}
{\it Proof.} We use a priori estimates on smooth solutions $u.$
Multiplying \eqref{e201} by $u$ and integrating over
$\mathbb{R}^{2},$ it is straightforward to see that the
$L^{2}(\mathbb{R}^{2})$-norm is conserved. Therefore, we only need
to show that
\begin{eqnarray}
\lefteqn{\partial_{t}\int \sum_{|j|\leq
N}[\,(\partial^{j}u_{xxx})^{2} +
(\partial^{j}(\partial_{x}^{-1}u_{yy}))^{2}\,]}\nonumber
\\
\label{e416}& \leq & c\left[\int \left(u^{2} +
\sum_{|j|\leq N}[\,(\partial^{j}u_{xxx})^{2} +
(\partial^{j}(\partial_{x}^{-1}u_{yy}))^{2}\,]\right) \right]^{3/2}.
\end{eqnarray}
We consider the case $j=(0,\,0).$ The case $j\neq (0,\,0)$ is
handled in a similar way. \\
\\
Applying $\partial_{x}^{3}$ to \eqref{e201} we obtain
\begin{eqnarray}
\label{e417}&  & u_{xxxt} + u_{xxxxxx} + u_{xxxx} +
(u\,u_{x})_{xxx} - u_{xxyy} = 0.
\end{eqnarray}
Multiplying \eqref{e417} by $2\,u_{xxx}$ and integrating over
$\mathbb{R}^{2}$ we obtain
\begin{eqnarray}
&  & 2\int u_{xxx}\,u_{xxxt} +
2\int u_{xxx}\,u_{xxxxxx} +
2\int u_{xxx}\,u_{xxxx} \nonumber \\
\label{e418}&  &
+\;2\int u_{xxx}\,(u\,u_{x})_{xxx} -
2\int u_{xxx}\,u_{xxyy} = 0.
\end{eqnarray}
Using in \eqref{e418} straightforward integration by parts, we
obtain
\begin{eqnarray*}
\lefteqn{\partial_{t}\int u_{xxx}^{2} = -\;2\int
u_{xxx}\,(u\,u_{x})_{xxx}} \\
& = & -\;2\int [\,3\,u_{xx}^{2} + 4\,u_{x}\,u_{xxx} +
u\,u_{xxxx}\,]\,u_{xxx}
 \\
& = & -\;7\int u_{x}\,u_{xxx}^{2}
 \leq 7\,||u_{x}||_{L^{\infty}(\mathbb{R}^{2})}
\int_{\mathbb{R}^{2}}u_{xxx}^{2} \\
& \leq & c\left(\int [\,u_{x}^{2} + u_{xxx}^{2} +
u_{xy}^{2}\,]\right)^{1/2}\int_{\mathbb{R}^{2}}u_{xxx}^{2}
\\
& \leq & c\left(\int [\,u^{2} + u_{xxx}^{2} +
(\partial_{x}^{-1}u_{yy}^{2})^{2}\,] \right)^{3/2}.
\end{eqnarray*}
In a similar way, but now apply
$\partial_{x}^{-1}\partial_{y}^{2}$ to \eqref{e201} instead of
$\partial_{x}^{3}$ and multiply by $2\,\partial_{x}^{-1}u_{yy}$
instead of $2\,u_{xxx}$ we get
\begin{eqnarray*}
\lefteqn{\partial_{t}\int (\partial_{x}^{-1}u_{yy})^{2}}\\
& \leq &
c\left|\int (u^{2})_{yy}\,(\partial_{x}^{-1}u_{yy})\right|
= c\left|\int [\,u_{yy}\,u +
u_{y}^{2}\,]\,(\partial_{x}^{-1}u_{yy})\right|\\
& \leq &
c\left|\int u_{x}\,(\partial_{x}^{-1}u_{yy})^{2}\right|
+
\left|\int u_{y}^{2}\,(\partial_{x}^{-1}u_{yy})\right|\\
& \leq & c||u_{x}||_{L^{\infty}(\mathbb{R}^{2})}
\int (\partial_{x}^{-1}u_{yy})^{2} +
\left(\int u_{y}^{4} \right)^{1/2}
\left(\int (\partial_{x}^{-1}u_{yy})^{2} \right)^{1/2}\\
& \leq & c||u_{x}||_{L^{\infty}(\mathbb{R}^{2})}
\int (\partial_{x}^{-1}u_{yy})^{2} +
\left(\int [\,u_{y}^{2} + u_{xy}^{2} +
(\partial_{x}^{-1}u_{yy})^{2}\,] \right)
\left(\int (\partial_{x}^{-1}u_{yy})^{2} \right)^{1/2}\\
& \leq & c\left(\int [\,u^{2} + u_{xxx}^{2} +
(\partial_{x}^{-1}u_{yy})^{2}\,]\right)^{3/2}.
\end{eqnarray*}
The lemma follows. \hfill $\square$ \\
\\
{\bf Corollary 4.4.} {\it Let $u$ be the solution to \eqref{e201}
with initial data $\phi\in Y^{N}(\mathbb{R}^{2}).$ Denote by
$0<T<+\infty$ the life span of this solution in
$Y^{N}(\mathbb{R}^{2}).$ Then there exists $0<T'\leq T,$ depending
only on the norm of $\phi\in X^{0}(\mathbb{R}^{2})$ such that
$u\in L^{\infty}([0,\,T']:\,X^{N}(\mathbb{R}^{2})).$}\\
\\
{\it Proof.} Let
\begin{eqnarray*}
h(t) = \int_{\mathbb{R}^{2}}\left(u^{2} + \sum_{|j|\leq
N}[\,(\partial^{j}u_{xxx})^{2} +
(\partial^{j}(\partial_{x}^{-1}u_{yy}))^{2}\,]\right)dx\,dy\equiv
||u||_{X^{N}}^2.
\end{eqnarray*}
Using \eqref{e415} we have $h'(t)\leq c\;[\,h(t)\,]^{3/2}.$
Integrating this inequality with respect to $t$, we obtain that
$h(t)^{1/2}\leq c/(h(0)^{-1/2} - t)$ and therefore, we get a lower
bound on the time of existence of $h(t)$ depending only on
$h(0).$ \hfill $\square$ \\
\\
{\bf Corollary 4.5.} {\it Let $\phi\in X^{N}(\mathbb{R}^{2})$ for
some $N\geq 0$ and let $\phi^{(n)}$ be a sequence converging to
$\phi$ in $X^{N}(\mathbb{R}^{2}).$ Let $u$ and $u^{(n)}$ be the
corresponding unique solutions, given by Theorems 4.1 and 4.2 and
Corollary 4.4, in $L^{\infty}([0,\,T]:\,X^{N}(\mathbb{R}^{2}))$ for
a time $T$ depending only on
$\sup_{n}||\phi^{(n)}||_{X^{0}(\mathbb{R}^{2})}.$ Then }
\begin{eqnarray}
\label{e419}u^{(n)}\stackrel{*}\rightharpoonup u\quad \mbox{weakly
in}\quad L^{\infty}([0,\,T]:\,X^{N}(\mathbb{R}^{2})).
\end{eqnarray}
{\it Proof.} By assumption $u\in
L^{\infty}([0,\,T]:\,X^{N}(\mathbb{R}^{2})),$ then there exists a
weak$^{*}$ convergent subsequence, still denoted $u^{(n)}$ such
that
\begin{eqnarray*}
u^{(n)}\stackrel{*}\rightharpoonup u\quad \mbox{weakly in}\quad
L^{\infty}([0,\,T]:\,X^{N}(\mathbb{R}^{2}))\hookrightarrow
L^{\infty}([0,\,T]:\,H^{1}(\mathbb{R}^{2})).
\end{eqnarray*}
Moreover, by equation \eqref{e201}, $u^{(n)}\in
L^{\infty}([0,\,T]:\,X^{N}(\mathbb{R}^{2}))$ implies
$u_{t}^{(n)}\in L^{\infty}([0,\,T]:\,L^{2}(\mathbb{R}^{2})).$ By
The Lions-Aubin compactness theorem,
\begin{eqnarray*}
u^{(n)}\rightarrow u\quad \mbox{strongly in}\quad
L^{\infty}([0,\,T]:\,H_{loc}^{1/2}(\mathbb{R}^{2})).
\end{eqnarray*}
Now we just need to show that each term in \eqref{e201} converges
to its correct limit, and $u_{t}^{(n)}\rightarrow u_{t}$ for $u\in
L^{\infty}([0,\,T]:\,X^{N}(\mathbb{R}^{2})).$\\
\\
The only thing we need to show is that the nonlinear term
converges to its correct limit, namely that
$u^{(n)}\,u_{x}^{(n)}\rightarrow u\,u_{x}.$ We know that
$u_{x}^{(n)}\stackrel{*}\rightharpoonup u_{x}$ weakly in $
L^{\infty}([0,\,T]:\,H^{1}(\mathbb{R}^{2}))$ and
$u^{(n)}\rightarrow u$ strongly in $
L^{\infty}([0,\,T]:\,H_{loc}^{1/2}(\mathbb{R}^{2})).$ Therefore,
their product converges in
$L^{2}([0,\,T]:\,L_{loc}^{1}(\mathbb{R}^{2})).$ Clearly, the
linear terms also converge in
$L^{2}([0,\,T]:\,L_{loc}^{1}(\mathbb{R}^{2}))$ and therefore, we
conclude that $u_{t}^{(n)}\rightarrow u_{t}$ in
$L^{2}([0,\,T]:\,L_{loc}^{1}(\mathbb{R}^{2})).$ The proof follows.
\hfill $\square$

\renewcommand{\theequation}{\thesection.\arabic{equation}}
\setcounter{equation}{0}\section{Estimate of error terms}
In this section we prove the main estimates used in our gain of regularity theorem.\\
\\
\nd{\bf Theorem 5.1}
{\it Let $L \geq 2$.  For $u$ a solution of (2.1), sufficiently smooth and with sufficient
decay at infinity,
\begin{equation}
  \sup_{0 \leq t \leq T} \int f_\alpha (\partial^\alpha u)^2 + \int_0^T \int g_\alpha
  (\partial^\alpha u_x)^2 \leq C
\end{equation}
for $L+1 \leq |\alpha| \leq 2L-1$, $2L - |\alpha|-\alpha_2 \geq 1$, where $f_\alpha
\in W_{\sigma,2L-|\alpha|-\alpha_2,|\alpha|-L}$, $g_\alpha \in W_{\sigma, 2L-|\alpha|
-\alpha_2-1,|\alpha|-L}$ and $C$ depends only on $||u||_{X^1}$ and
\begin{align}
\label{allowed-term-1}
& \sup_{0 \leq t \leq T} \int f_\gamma(\partial^\gamma u)^2 \\
\label{allowed-term-2}
& \int_0^T \int g_\gamma (\partial^\gamma u_x)^2
\end{align}
where $\gamma = (\gamma_1,\gamma_2) \in \mathbb Z^+ \times \mathbb Z^+$, $|\gamma|
\leq |\alpha|-1$, $f_\gamma \in
W_{\sigma, 2L - |\gamma|-\gamma_2,|\gamma|-L}, g_\gamma \in W_{\sigma,2L-|\gamma|-
\gamma_2-1,|\gamma|-L}$ for $|\gamma| \geq L$, $2L-|\gamma|-\gamma_2 \geq 1$ and
$f_\gamma \in W_{0, \gamma_1,0},
g_\gamma \in W_{\sigma, \gamma_1-1,0}$ for $0 \leq |\gamma| \leq L$.}\\
\\
\\
The idea of the proof is the following.  For a given $\alpha$ satisfying the hypotheses above,
we choose a weight function $f_\alpha \approx t^{|\alpha|-L}
x^{2L-|\alpha|-\alpha_2}$ for $x > 1$ and $f_\alpha \approx t^{|\alpha|-L}e^{\sigma x}$
for $x < -1$. Then with this choice of weight function,
we apply the operator $\partial^\alpha$ to
\eqref{e201}, multiply by $f_\alpha\pa u$ and integrate over $\mathbb R^2$ to obtain the
main equality stated in \eqref{e203}.  In this theorem, we bound the last three terms on the
left-hand side of \eqref{e203} by terms of the form \eqref{allowed-term-1} and \eqref{allowed-term-2}.  \\

\nd{\it Proof. }  For each $\alpha$ we apply the operator $\partial^\alpha$ to (2.1),
multiply our differentiated equation by $2 \fa (\pa u)$ where we take
\begin{equation}
  \fa(x,t) = \int_{-\infty}^x g_\alpha(z,t) \,dz \quad \text{for } g_\alpha \in
  W_{\sigma, 2L-|\alpha|-\alpha_2-1,|\alpha|-L},
\end{equation}
and integrate over $\mathbb R^2 \times [0,t]$ for $0 \leq t \leq T$.
As stated in Lemma 2.1, we arrive at our main equality
\begin{align*}
& \partial_t \int \fa (\pa u)^2 + 3 \int (\fa)_x (\pa u_x)^2 - \int [(\fa)_t + (\fa)_{xxx} +
(\fa)_x] (\pa u)^2 \\
& \qquad - \int (\fa)_x (\pa \partial_x^{-1}u_y)^2 + 2 \int \fa (\pa u)\pa(uu_x) = 0.
\end{align*}
Using (1.5) and $f(\cdot,0) = 0$ we get the following identity after integrating with respect
to $t$,
\begin{equation}
\label{before-estimating}
\begin{split}
& \int \fa (\pa u)^2 + 3 \int_0^T \int (\fa)_x (\pa u_x)^2 \\
& \leq \int_0^T \int (\fa)_x (\pa \partial_x^{-1} u_y)^2 + C \int_0^T \int \fa (\pa u)^2
 + 2 \left|\int_0^T \int \fa (\pa u) \pa (uu_x) \right|.
\end{split}
\end{equation}
We notice that the first term on the right-hand side of \eqref{before-estimating} can be
written as
\begin{equation}
\label{inverse-term}
\int_0^T \int (\fa)_x (\pa \partial_x^{-1} u_y)^2 = C \int_0^T \int t g_\gamma (\partial^\gamma u_x)^2
\end{equation}
for some $g_\gamma \in W_{\sigma,2L-|\gamma|-\gamma_2-1,|\gamma|-L}$
where $\gamma = (\alpha_1-2,\alpha_2+1)$.  Further, we notice that $2L-|\gamma|-\gamma_2 \geq 1$
and $\alpha_1 \geq 2$ since $2L-|\alpha|-\alpha_2\geq 1$ and $L+1 \leq |\alpha|$.
Therefore, \eqref{inverse-term} is of the form specified by \eqref{allowed-term-2}.
Therefore,
\begin{equation}
\label{after-inverse-term}
\begin{split}
& \int \fa (\pa u)^2 + 3 \int_0^T \int (\fa)_x (\pa u_x)^2 \\
& \qquad \leq C + C \int_0^T \int \fa (\pa u)^2
 + 2 \left|\int_0^T \int \fa (\pa u) \pa (uu_x) \right|
\end{split}
\end{equation}
where $C$ depends only on terms of the form \eqref{allowed-term-2}.
We now need to estimate the term
\[
\left|\int_0^T \int \fa (\pa u)\pa (uu_x)\right|.
\]
Each term is of the form
\[
\left|\int_0^T \int \fa (\pa u)(\partial^r u)(\partial^s u_x)\right|.
\]
where $r_1 + s_1 = \alpha_1$, $r_2 + s_2 = \alpha_2$.  Below we consider all terms of level
$|\alpha|$.

\vskip 5pt
\nd{\bf The case $|s| = |\alpha|$.}  In this case, $|r|=0$, and we have
\begin{align*}
  \left|\int_0^T \int \fa (\pa u)(\partial^r u)(\partial^s u_x)\right| & = \left|\int_0^T
  \int \fa (\pa u) u (\pa u_x)\right| \\
  & = C \left|\int_0^T \int (\fa u)_x (\pa u)^2\right| \\
  & \leq C ||u||_{X^0} \int_0^T \int \fa (\pa u)^2.
\end{align*}

\vskip 5pt
\nd{\bf The case $|s| = |\alpha|-1$.}  In this case, $|r|=1$ giving us the following two subcases:

(a) {\bf The subcase $r = (1,0)$.}  In this case, we have
\begin{align*}
  \left|\int_0^T \int \fa (\pa u)(\partial^r u)(\partial^s u_x)\right| & = \left|\int_0^T \int
  \fa (\pa u)u_x (\pa u)\right| \\
  & \leq C||u||_{X^0} \int_0^T \int \fa (\pa u)^2.
\end{align*}

(b) {\bf The subcase $r = (0,1)$.}  We note that this case will only occur if
$\alpha_2 \geq 1$.  In this case, we have
\begin{align*}
  \left|\int_0^T \int \fa (\pa u)(\partial^r u)(\partial^s u_x)\right| & = \left|\int_0^T
  \int \fa (\pa u) u_y (\partial_x^{\alpha_1} \partial_y^{\alpha_2-1}u_x)\right| \\
  & \leq C||u||_{X^1} \left(\int_0^T \int \fa (\pa u)^2\right)^{1/2} \left(\int_0^T
  \int \fa (\partial_x^{\alpha_1+1}\partial_y^{\alpha_2-1}u)^2\right)^{1/2}
\end{align*}
Since $f_\alpha \approx x^{2L-|\alpha|-\alpha_2}$ as $x \rightarrow \infty$, it is clear that
$f_\alpha \leq C f_{\alpha_1+1,\alpha_2-1}$

\nd{\bf The case $|s| = |\alpha|-2$.}  We have three subcases to consider.

(a) {\bf The subcase $r = (2,0)$.}  In this case, we have
\begin{align*}
  \left|\int_0^T \int \fa (\pa u)(\partial^r u)(\partial^s u_x)\right| & = \left|\int_0^T \int
  \fa (\pa u)u_{xx} (\partial_x^{\alpha_1-2}\partial_y^{\alpha_2} u_x)\right| \\
  & \leq C||u_{xx}||_{L^\infty} \left(\int_0^T \int \fa (\pa u)^2\right)^{1/2} \left(\int_0^T
  \int \fa (\partial_x^{\alpha_1-1}\partial_y^{\alpha_2}u)^2\right)^{1/2} \\
  & \leq C ||u||_{X^1}\left(\int_0^T \int \fa (\pa u)^2\right)^{1/2} \left(\int_0^T
  \int \fa (\partial_x^{\alpha_1-1}\partial_y^{\alpha_2}u)^2\right)^{1/2}.
\end{align*}
The last term on the right-hand side above is of order $|\alpha|-1 \geq L$.
The weight function $f_\alpha \in W_{\sigma, 2L-|\alpha|-\alpha_2,|\alpha|-L}$.
Since $2L-|\alpha|-\alpha_2 < 2L-(\alpha_1-1+\alpha_2)
-\alpha_2$, we see that this term is bounded by a term of the form \eqref{allowed-term-1}.

(b) {\bf The subcase $r=(1,1)$.}  In this case, we have
\begin{align*}
  \left|\int_0^T \int \fa (\pa u)(\partial^r u)(\partial^s u_x)\right| & = \left|\int_0^T
  \int \fa (\pa u)u_{xy} (\partial_x^{\alpha_1-1}\partial_y^{\alpha_2-1}u_x)\right| \\
  & \leq C||u_{xy}||_{L^\infty} \left(\int_0^T \int \fa (\pa u)^2\right)^{1/2} \left(\int_0^T
  \int \fa (\partial_x^{\alpha_1}\partial_y^{\alpha_2-1} u)^2\right)^{1/2} \\
  & \leq C||u||_{X^1} \left(\int_0^T \int \fa (\pa u)^2\right)^{1/2} \left(\int_0^T
  \int \fa (\partial_x^{\alpha_1}\partial_y^{\alpha_2-1} u)^2\right)^{1/2}
\end{align*}
Using the fact that $\fa \in W_{\sigma,2L-|\alpha|-\alpha_2,|\alpha|-L}$ and $2L-|\alpha|
-\alpha_2 < 2L-(\alpha_1+\alpha_2-1) - (\alpha_2-1)$, we conclude that the last term
is bounded by a term of the form \eqref{allowed-term-1}.

(c) {\bf The subcase $r=(0,2)$.}  In this case, we have
\begin{align*}
  \left|\int_0^T \int \fa (\pa u)(\partial^r u)(\partial^s u_x)\right| & = \left|\int_0^T
  \int \fa (\pa u)u_{yy} (\partial_x^{\alpha_1} \partial_y^{\alpha_2-2}u_x)\right| \\
  & \leq \left(\int_0^T \int \fa (\pa u)^2\right)^{1/2} \left(\int_0^T \int u_{yy}^4
  \right)^{1/4} \\
  & \qquad \times \left(\int_0^T \int (\fa^{1/2} \partial_x^{\alpha_1+1}\partial_y^{\alpha_2-2}u)^4
  \right)^{1/4}.
\end{align*}
Now
\[
\left(\int_0^T \int u_{yy}^4\right)^{1/4} \leq C \left(\int_0^T \int [u_{yy}^2 + u_{xyy}^2 +
(\partial_x^{-1} u_{yyy})^2]\right)^{1/2} \leq C||u||_{X^1}.
\]
Further,
\begin{align*}
  & \left(\int_0^T \int (\fa^{1/2} \partial_x^{\alpha_1+1}\partial_y^{\alpha_2-2}u)^4 \right)^{1/4}\\
  & \leq C \left(\int_0^T \int \fa (\partial_x^{\alpha_1+1} \partial_y^{\alpha_2-2}u)^2
  + \fa (\partial_x^{\alpha_1+2}\partial_y^{\alpha_2-2}u)^2 + \fa (\partial_x^{\alpha_1}
  \partial_y^{\alpha_2-1}u)^2    \right)^{1/2}
\end{align*}
Now the first and third terms in the integrand are of order $|\alpha|-1$ and are clearly bounded
by terms of the form \eqref{allowed-term-1}.  The second term in the integrand is
of order $|\alpha|$.  It will be bounded using Gronwall's inequality (and using the fact
that the order of the $y$ derivative is less than $\alpha_2$ and therefore this terms can handle
an even greater power of $x$.)

\vskip 5pt
\nd{\bf The case $|s| = |\alpha|-3$. }  In this case $|r|=3$.  First, we consider the case in
which $L \geq 4$.  Since $|\alpha| \geq L+1$, we note that $|s|+2 = |\alpha|-1 \geq L$.
Using this fact, we bound as follows.
\begin{align*}
  \left|\int_0^T \int \fa (\pa u)(\partial^r u)(\partial^s u_x)\right| & \leq \left(\int_0^T
  \int \fa (\pa u)^2\right)^{1/2} \left(\int_0^T \int (\partial^r u)^4\right)^{1/4} \\
  & \qquad \times 
  \left(\int_0^T \int (\fa^{1/2} \partial^s u_x)^4\right)^{1/4}
\end{align*}
Now
\begin{align*}
  \left(\int_0^T \int (\partial^r u)^4\right)^{1/4} \leq C\left(\int_0^T \int
  [(\partial^r u)^2 + (\partial^r u_x)^2 + (\partial^r u_y)^2]\right)^{1/2}.
\end{align*}
Since $|r| = 3$, each of these terms is at most of order $4 \leq L \leq |\alpha|-1$,
and, therefore, bounded by terms of the form \eqref{allowed-term-1}.
Similarly,
\begin{equation}
\label{s-alpha-3}
  \left(\int_0^T \int (\fa^{1/2} \partial^s u_x)^4\right)^{1/4} \leq C \left(\int_0^T \int
  \fa [(\partial^s u_x)^2 + (\partial^s u_{xx})^2 + (\partial^s u_{xy})^2]\right)^{1/2}
\end{equation}
We notice that the last two terms on the right-hand side of \eqref{s-alpha-3} are
of order $|s|+2=|\alpha|-1\geq L$.
In order to verify that we have the correct power of $x$, we note that
\[
2L-|\alpha|-\alpha_2 \leq 2L-(|s|+2)-(s_2+1)
\]
since $|s|= |\alpha|-3$.  Therefore, we can conclude that each of those
terms is bounded by a term of the form \eqref{allowed-term-1}.  Finally, we look at the first
term on the right-hand side of \eqref{s-alpha-3}.  If $|s|+1 \geq L$, then this term is bounded
by \eqref{allowed-term-1} as the other two terms.  If $|s|+1 < L$, then using the fact that
$|s|+2 \geq L$, we conclude that $|s| = L-2$, and, therefore,
\begin{align*}
2L-|\alpha|-\alpha_2 & \leq 2L-(|s|+2)-(s_2+1) \\
&\leq s_1+1
\end{align*}
Therefore, we conclude that the first term above is bounded by a term of the form
\eqref{allowed-term-1} of order $|s|+1 < L$.

We now look at the cases when $L=2$ or $L=3$.  In either case, if $|\alpha| \geq 5$, then
we can handle as above.  We first consider the case when $|\alpha|=4$ ($L=2$ or $L=3$).
In this case, using the fact that $|r|=3$ and $|s|=1$, we have
\begin{align*}
  \left|\int_0^T \int \fa (\pa u)(\partial^r u)(\partial^s u_x)\right| & \leq
  C||\partial^s u_x||_{L^\infty} \left(\int_0^T
  \int \fa (\pa u)^2\right)^{1/2} \left(\int_0^T \int \fa (\partial^r u)^2 \right)^{1/2} \\
  & \leq C||u||_{X^1} \left(\int_0^T \int \fa (\pa u)^2\right)^{1/2} \left(\int_0^T
  \int \fa (\partial^r u)^2\right)^{1/2}.
\end{align*}
Since $|r|=3=|\alpha|-1$ and $2L-|\alpha|-\alpha_2 \leq 2L-|r|-r_2$, we see that the last
term above is bounded by a term of the form \eqref{allowed-term-1}.

Last, we consider $|\alpha|=3$.  In this case, we must have $L=2$, $|r|=3$ and $|s|=0$.
Therefore, $r = \alpha$.  We bound as follows:
\begin{align*}
  \left|\int_0^T \int \fa (\pa u)(\partial^r u)(\partial^s u_x) \right| & \leq C||u_x||_{L^\infty}
  \int_0^T \int \fa (\pa u)^2 \\
  & \leq C||u||_{X^1} \int_0^T \int \fa (\pa u)^2.
\end{align*}

\nd{\bf The case $|s| \leq |\alpha|-4$.}  We consider the set $A = \{x: x > 1\}$.
The set $A_{-1} = \{x < -1\}$ can be handled similarly.  We have
\begin{align*}
  & \left|\int_0^T \int_A \fa (\pa u)(\partial^r u)(\partial^s u_x)\right| \\
  & \leq C T^M||
  t^{\nu_s}\partial^s u_x||_{L^\infty(A)} \left(\int_0^T \int_A \fa (\pa u)^2\right)^{1/2}
  \left(\int_0^T \int_A t^{\nu_r} x^{2L-|\alpha|-\alpha_2}(\partial^r u)^2\right)^{1/2}
\end{align*}
where $\nu_s = \frac{(|s|+3-L)^+}2$ and $\nu_r = \frac{(|r|-L)^+}2$.  First, we must verify that
$M \geq 0$.  We see that
\begin{align*}
M & = \frac {|\alpha|-L}2 - \nu_s - \nu_r \\
& = \frac{|\alpha|-L-(|s|+3-L)^+-(|r|-L)^+}2 \\
& \geq \frac {|\alpha|-|s|-|r|-3+L}2 \\
& = \frac {L-3}2 \geq 0
\end{align*}
as long as $L \geq 3$.  By assumption, $L \geq 2$.  If $L = 2$, then $|\alpha| = 3$.  In
that case, we cannot have $|s| \leq |\alpha| - 4$.  Therefore, we conclude that $M \geq 0$.
Further,
\begin{align*}
  ||t^{\nu_s}\partial^s u_x||_{L^\infty} & \leq C \left(\int t^{2\nu_s} [(\partial^su_x)^2
  + (\partial^s u_{xxx})^2 + (\partial^s u_{xy})^2]\right)^{1/2}.
\end{align*}
Each of those terms is of order at most $|\alpha|-1$, and therefore,
bounded by terms of the form \eqref{allowed-term-1}.  Further,
$2L-|\alpha|-\alpha_2 \leq 2L-|r|-r_2$.  Therefore, for $|r| \leq
|\alpha|-1$, the last term above is bounded by terms of the form
\eqref{allowed-term-1}. In the case that $|r|=|\alpha|$, we have
$|s| = 0$, and therefore
\begin{align*}
  \left|\int_0^T \int \fa (\pa u)(\partial^r u)(\partial^s u_x)\right| & \leq C ||u_x||
  _{L^\infty} \int_0^T \int \fa (\pa u)^2 \\
  & \leq C ||u||_{X^0} \int_0^T \int \fa (\pa u)^2.
\end{align*}

Combining our estimates above on
\[
\left|\int_0^T \int \fa (\pa u)(\partial^r u)(\partial^s u_x)\right|
\]
with \eqref{after-inverse-term}, we see that
\begin{equation}
\begin{split}
\label{power-1} \int \fa (\pa u)^2 + 3 \int_0^T \int (\fa)_x (\pa
u_x)^2 & \leq C + C\sum_{\stackrel{|\gamma| = |\alpha|}{
2L-|\gamma|-\gamma_2 \geq 1}} \int_0^T \int f_\gamma
(\partial^\gamma u)^2
\end{split}
\end{equation}
where the constant $C$ depends only on terms of the form \eqref{allowed-term-1} and
\eqref{allowed-term-2}.

Using the above estimate for all derivatives $\gamma$ of order $|\alpha|$ such that
$2L-|\gamma|-\gamma_2 \geq 1$, we see that
\begin{equation}
\begin{split}
\sum_{\stackrel{|\gamma|=|\alpha|}{2L-|\gamma|-\gamma_2 \geq 1}}&
\left[\int f_\gamma (\pa u)^2 +
3 \int_0^T \int (f_\gamma)_x (\pa u_x)^2\right] \\
& \leq C + \sum_{\stackrel{|\gamma|=|\alpha|}{2L-|\gamma|-\gamma_2
\geq 1}}\int_0^T \int f_\gamma (\partial^\gamma u)^2.
\end{split}
\end{equation}
where $C$ depends only on terms of the form \eqref{allowed-term-1} and \eqref{allowed-term-2}.
Applying Gronwall's inequality, we get the desired estimate. \hfill $\square$ \\
\\

\renewcommand{\theequation}{\thesection.\arabic{equation}}
\setcounter{equation}{0}\section{Persistence Theorem}
In section four we proved the existence of a solution $u$ to \eqref{e201} in $L^\infty([0,T];
X^N(\mathbb R^2))$ for given initial data $\phi \in X^N(\mathbb R^2)$.  In this section,
we prove that if, in addition, our initial data $\phi$ lies in the weighted space
$\widetilde H_x^K(W_{0\;K\;0})$ for some $K \geq 0$, then the solution
$u$ also lies in $L^\infty([0,T];\widetilde H_x^K(W_{0\;K\;0}))$.  This property
is known as a ``persistence" property of the initial data.  This property provides a basis
for starting the induction in our Gain of Regularity theorem in Section 7.\\
\\
\nd{\bf Theorem 6.1}
{\it Suppose $u \in L^\infty([0,T]:X^1(\mathbb R^2))$ with initial data $\phi(x,y)
\in X^1(\mathbb R^2)$ such that $\phi$ also lies in $\widetilde H_x^K(W_{0\;K\;0})$
for some integer $K \geq 0$.  Then
\[
u \in L^\infty([0,T]:X^1(\mathbb R^2) \cap \widetilde H_x^K(W_{0\;K\;0}))
\]
and
\[
\sup_{0 \leq t \leq T}\int \faone(\pa u)^2 + \int_0^T \int \gaone (\pa u_x)^2 \,dt \leq C
\]
for $|\alpha| \leq K$, $\alpha_1 \neq 0$, where $\faone \in W_{0\;\alpha_1\;0}$ and
$\gaone \in W_{\sigma\; \alpha_1-1\; 0}$ for $\sigma > 0$
arbitrary and $C$ depends only on $T$ and the norm of $\phi \in X^1(\mathbb R^2)
\cap \widetilde H_x^K(W_{0\;K\;0})$.}


\vskip 5pt
\nd{\bf Proof. }
We use induction on $j=|\alpha|$ for $1 \leq j \leq K$.  The case that $j = 0$ follows
from conservation of $L^2$ norm.
We derive formally some a priori estimates for the solution where the bound involves only
the norms of $u \in L^\infty([0,T]:X^1(\mathbb R^2))$ and the norms of $\phi \in \widetilde H_x^K(
W_{0\;K\;0})$.  Then, we can apply convergence arguments
to show that the result holds true for general solutions.  In order to do so, we need to
approximate general solutions $u \in X^1(\mathbb R^2)$ by smooth solutions and approximate
general weight functions $f \in W_{0\;j\;0}$ by smooth, bounded weight functions.  The first
of these procedures has already been discussed, so we will concentrate on the second.

For a fixed $i$, we begin by taking a sequence of bounded weight functions $g_{i,\delta}$
which decay as $|x| \rightarrow \infty$ and which approximate $g_i \in W_{\sigma\; i-1\; 0}$
with $\sigma > 0$ from below, uniformly on any half-line $(-\infty,c)$.  Define
the weight functions
\[
f_{i,\delta}(x,t) = 1 + \int_{-\infty}^x g_{i,\delta}(z,t) \,dz.
\]
Therefore, the functions $f_{i,\delta}$ are bounded weight functions approximating
$f_i \in W_{0,i,0}$ from below, uniformly on compact sets.

>From (5.3) and using the fact that $\partial_t (\fid) \leq c \fid$ and $\partial_x (\fid) \leq c
\fid$, we have
\begin{align*}
  & \int \fid (\pa u)^2 + 3 \int_0^T \int (\fid)_x (\pa u_x)^2 \,dt
   \leq \int_0^T (\fid)_x (\pa \partial_x^{-1} u_y)^2 \,dt \\
   & \qquad + C\int_0^T
  \int \fid (\pa u)^2 \,dt + 2 \left|\int_0^T \int \fid (\pa u) \pa (uu_x) \,dt \right|.
\end{align*}

\vskip 5pt
\nd {\bf The case $j=1$.}
\vskip 5pt

(a) {\bf The subcase $\alpha = (1,0)$.}
Defining $g_{1,\delta}$ and $f_{1,\delta}$ as above,
we see that $f_{1,\delta}$ will approximate $f_1 \in W_{0\;1\;0}$ from below.
Differentiating (2.1) in the $x-$variable, multiplying by $2f_{1,\delta}$ and integrating
over $\mathbb R^2$, we have
\begin{equation}
\begin{split}
\label{alpha1beta0}
& \partial_t \int \foned u_x^2 + 3 \int(\foned)_xu_{xx}^2 \\
& \qquad = \int (\foned)_x u_y^2 + \int [\partial_t\foned + \partial_x^3\foned + \partial_x
\foned]u_x^2 - 2 \int \foned u_x(uu_x)_x \\
& \leq C \int u_y^2 + C \int \foned u_x^2 + 2 \left|\int \foned u_x(uu_x)_x\right|
\end{split}
\end{equation}

Moreover
\begin{align*}
  \left|\int \foned u_x(uu_x)_x \right|
  & = \left|\int \foned u_x[u_x^2 + uu_{xx}]\right| \\
  & = \left|\int \foned [u_x^3 + uu_{x}u_{xx}] \right| \\
  & \leq c ||u_x||_{L^\infty} \int \foned u_x^2 + \frac 12\left|\int (\foned u)_x u_x^2 \right| \\
  & \leq c(||u||_{L^\infty} + ||u_x||_{L^\infty})\int \foned u_x^2 \\
  & \leq C ||u||_{X^0}\int \foned u_x^2.
\end{align*}
Combining this estimate with \eqref{alpha1beta0}, we conclude that for $0 \leq t \leq T$,
\[
\int \foned(\cdot,t) u_x^2 + 3 \int_0^T \int (\foned)_x u_{xx}^2 \leq C \int \foned(\cdot,0)
\phi_x^2 + C \int_0^T \int u_y^2 + C \int_0^T \int \foned u_x^2.
\]
Applying Gronwall's inequality, we conclude that
\[
\sup_{0 \leq t \leq T} \int \foned u_x^2 + 3 \int_0^T \int (\foned)_x u_{xx}^2 \leq C
\]
where $C$ does not depend on $\delta$ but only on $T$ and the norm of $\phi \in
X^1(\mathbb R^2) \cap \widetilde H_x^1(W_{0\;1\;0})$.
Taking the limit as $\delta\rightarrow \infty$, we conclude that
\begin{equation}
  \sup_{0 \leq t \leq T} \int f_1 u_x^2 + 3 \int_0^T \int g_1 u_{xx}^2 \leq C,
\end{equation}
as claimed.

\vskip 5pt
 (b) {\bf The subcase $\alpha = (0,1)$.}
Here, our weight function $f_0 \in W_{0\;0\;0}$.
Differentiating (2.1) in the $y$-variable, multiplying by $2\fzerod u_y$ and
integrating over $\mathbb R^2$, we have
\[
2 \int \fzerod u_y u_{yt} + 2 \int \fzerod u_y u_{xxxy} - 2 \int \fzerod u_y \partial_x^{-1}
u_{yyy} + 2 \int \fzerod u_y u_{xy} + 2 \int \fzerod u_y (uu_x)_y = 0.
\]
Integrating each term by parts gets
\begin{equation}
  \begin{split}
    \label{alpha0beta1}
    & \partial_t \int \fzerod u_y^2 + 3 \int (\fzerod)_x u_{xy}^2 \\
    & \qquad \leq \int (\fzerod)_x (\partial_x^{-1}u_{yy})^2 + C \int \fzerod u_y^2
    + 2 \left|\int \fzerod u_y (uu_x)_y\right| \\
    & \leq C \int (\partial_x^{-1} u_{yy})^2 + C \int \fzerod u_y^2 + 2 \left|
    \int \fzerod u_y (uu_x)_y\right|.
  \end{split}
\end{equation}
Moreover,
\begin{align*}
  \left|\int \fzerod u_y (uu_x)_y\right| & = \left|\int \fzerod u_y (u_yu_x + uu_{xy})\right| \\
  & \leq C ||u_x||_{L^\infty} \int \fzerod u_y^2 + C ||u||_{L^\infty} \int (\fzerod)_x u_y^2 \\
  & \leq C ||u||_{X^0}\int \fzerod u_y^2.
\end{align*}
Combining this estimate with \eqref{alpha0beta1}, we conclude that for $0 \leq t \leq T$,
\begin{equation}
  \int \fzerod(\cdot,t)u_y^2 + 3 \int_0^T \int (\fzerod)_x u_{xy}^2 \leq \int \fzerod(\cdot,0)
  \phi_y^2 + C \int_0^T \int (\partial_x^{-1}u_{yy})^2 + C \int_0^T \int \fzerod u_y^2.
\end{equation}
Applying Gronwall's inequality, we conclude that
\[
\sup_{0 \leq t \leq T} \int \fzerod u_y^2 + 3 \int_0^T \int (\fzerod)_x u_{xy}^2 \leq C
\]
where $C$ does not depend on $\delta$, but only on $T$ and the norm of $\phi \in X^0(\mathbb R^2)
\cap \widetilde H_x^1(W_{0\;1\;0})$.  Passing to the limit, we conclude that
\begin{equation}
  \sup_{0 \leq t \leq T} \int f_0 u_y^2 + 3 \int_0^T \int g_0 u_{xy}^2 \leq C.
\end{equation}

\vskip 5pt

\nd{\bf The case $j = 2$.}

\vskip 5pt
(a) {\bf The subcase $\alpha = (2,0)$.}
In this case, $\ftwod$ will approximate $f_2 \in W_{0\;2\;0}$.  In a similar way as above, we have
\begin{align*}
  & 2 \int \ftwod u_{xx} u_{xxt} + 2 \int \ftwod u_{xx}u_{xxxxx} - 2 \int \ftwod u_{xx} u_{xyy} \\
  & \qquad + 2 \int \ftwod u_{xx}u_{xxx} + 2 \int \ftwod u_{xx}(uu_x)_{xx} = 0.
\end{align*}
Integrating each term by parts gets
\begin{align*}
  & \partial_t \int \ftwod u_{xx}^2 + 3 \int (\ftwod)_xu_{xxx}^2 \\
  & \qquad = c \int (\ftwod)_x u_{xy}^2 + \int [\partial_t \ftwod + \partial_x^3 \ftwod
  + \partial_x \ftwod]u_{xx}^2 - 2 \int \ftwod u_{xx} (uu_x)_{xx} \\
  & \leq \int \ftwod u_{xy}^2 + c \int \ftwod u_{xx}^2 + 2 \left|\int \ftwod u_{xx}
  (uu_x)_{xx} \right|.
\end{align*}
Moreover
\begin{align*}
  \left|\int \ftwod u_{xx}(uu_x)_{xx} \right|
  & = \left|\int \ftwod u_{xx}[3u_{x} u_{xx} + uu_{xxx}] \right| \\
  & = \left|\int \ftwod[3u_{x} u_{xx}^2 + uu_{xx} u_{xxx}] \right| \\
  & \leq c||u_x||_{L^\infty} \int \ftwod u_{xx}^2 + \frac{1}{2}\left|(\ftwod u)_x u_{xx}^2 \right| \\
  & \leq c(||u||_{L^\infty} + ||u_x||_{L^\infty})\int \ftwod u_{xx}^2 \\
  & \leq C ||u||_{X^0} \int \ftwod u_{xx}^2.
\end{align*}
Therefore,
\[
\partial_t \int \ftwod u_{xx}^2 + 3 \int (\ftwod)_x u_{xxx}^2 \leq \int (\ftwod)_x u_{xy}^2
+ C \int \ftwod u_{xx}^2,
\]
where $C$ depends only on the norm of $\phi \in X^0(\mathbb R^2)$.
We will combine this estimate with the estimate below.

\vskip 5pt
 (b) {\bf The subcase $\alpha = (1,1)$}
Applying $\partial_x\partial_y$ to (2.1), multiplying by $2 \foned u_{xy}$ where
$\foned$ approximates $f_1 \in W_{0\;1\;0}$, and integrating over $\mathbb R^2$, we have
\begin{align*}
& 2 \int \foned u_{xy} u_{xyt} + 2 \int \foned u_{xy}u_{xxxxy} - 2 \int \foned u_{xy}u_{yyy} \\
& \qquad + 2 \int \foned u_{xy}u_{xxy} + 2 \int \foned u_{xy}(uu_x)_{xy} = 0.
\end{align*}
Integrating each term by parts gets
\begin{align*}
  & \partial_t \int \foned u_{xy}^2 + 3 \int (\foned)_x u_{xxy}^2 \\
  & \qquad = \int (\foned)_x u_{yy}^2 + \int [\partial_t \foned + \partial_x^3 \foned +
  \partial_x \foned]u_{xy}^2 - 2 \int \foned u_{xy}(uu_x)_{xy} \\
  & \leq \int (\foned)_x u_{yy}^2 + c \int \foned u_{xy}^2 + 2 \left|\int \foned u_{xy}
  (uu_x)_{xy}\right|.
\end{align*}
Moreover,
\begin{align*}
  \left|\int \foned u_{xy}(uu_x)_{xy}\right|
  & = \left|\int \foned u_{xy}(2u_x u_{xy} + u_{xx}u_y + uu_{xxy})\right| \\
  & \leq C(||u_x||_{L^\infty} + ||u||_{L^\infty})\int \foned u_{xy}^2 +
  \left|\int \foned u_{xx}u_{xy} u_y\right| \\
  & \leq C||u||_{X^0}\int \foned u_{xy}^2 + C||u_y||_{L^\infty} (\int \foned u_{xx}^2
  + \int \foned u_{xy}^2) \\
  & \leq C||u||_{X^0} \int \foned u_{xy}^2 + C ||u||_{X^1} \int \foned (u_{xx}^2 + u_{xy}^2)
\end{align*}
since
\begin{align*}
||u_y||_{L^\infty} & \leq \left(\int u_y^2 + u_{xxy}^2 + u_{yy}^2\right)^{1/2} \leq
||u||_{X^1}.
\end{align*}
Therefore,
\[
\partial_t \int \foned u_{xy}^2 + 3 \int (\foned)_x u_{xxy}^2 \leq \int (\foned)_x u_{yy}^2
+ C \int \foned u_{xx}^2 + C \int \foned u_{xy}^2.
\]

\vskip 5pt
 (c) {\bf The subcase $\alpha = (0,2)$.}
Applying $\partial_y^2$ to (2.1), multiplying by
$u_{yy}$ and integrating over $\mathbb R^2$, we have
\begin{align*}
& 2 \int u_{yy} u_{yyt} + 2 \int u_{yy}u_{xxxyy} - 2 \int u_{yy} \partial_x^{-1}u_{yyyy} \\
& \qquad + 2 \int u_{yy} u_{xyy} + 2 \int u_{yy}(uu_x)_{yy} = 0.
\end{align*}
Integrating by parts gets
\[
\partial_t \int u_{yy}^2 \leq 2 \left|\int u_{yy}(uu_x)_{yy}\right|.
\]
Now
\begin{align*}
  \left|\int u_{yy} (uu_x)_{yy}\right| & = \left|\int u_{yy}(u_{yy}u_x + 2u_y u_{xy} + uu_{xyy})
  \right| \\
  & \leq ||u_x||_{L^\infty} \int u_{yy}^2 + ||u_y||_{L^\infty} \left(\int u_{yy}^2 + u_{xy}^2
  \right) \\
  & \leq ||u||_{X^0}\int u_{yy}^2 + C||u||_{X^1} \int (u_{yy}^2 + u_{xy}^2).
\end{align*}

Now combining these estimates from (a), (b) and (c) above, we have
\begin{align*}
  & \partial_t \int (\ftwod u_{xx}^2 + \foned u_{xy}^2 + u_{yy}^2) + 3 \int \left[(\ftwod)_x
  u_{xxx}^2 + (\foned)_x u_{xxy}^2\right] \\
  & \qquad \leq C \int (\ftwod + \foned) u_{xx}^2
  + C \int ((\ftwod)_x + \foned + 1) u_{xy}^2 + \int ((\foned)_x + 1) u_{yy}^2,
\end{align*}
where $C$ depends only on the norm of $\phi \in X^1(\mathbb R^2)$.
Since $\ftwod$ approximates $f_2 \in W_{0\;2\;0}$ and $\foned$ approximates $f_1
\in W_{0\;1\;0}$, we can choose $\ftwod, \foned$ such that $(\ftwod)_x \leq C \foned$, etc.
Therefore,
\begin{align*}
  & \partial_t \int (\ftwod u_{xx}^2 + \foned u_{xy}^2 + u_{yy}^2) + 3 \int \left[(\ftwod)_x
  u_{xxx}^2 + (\foned)_x u_{xxy}^2\right] \\
  & \qquad \leq C \int (\ftwod u_{xx}^2+ \foned u_{xy}^2 + u_{yy}^2).
\end{align*}
Integrating with respect to $t$, we have
\begin{align*}
  & \int (\ftwod(\cdot,t) u_{xx}^2 + \foned(\cdot,t) u_{xy}^2 + u_{yy}^2) +
  3 \int_0^t\int \left[(\ftwod)_x u_{xxx}^2 + (\foned)_x u_{xxy}^2\right] \\
  & \leq \int (\ftwod(\cdot,0) \phi_{xx}^2 + \foned(\cdot,0) \phi_{xy}^2 + \phi_{yy}^2)
  + C \int_0^t \int (\ftwod u_{xx}^2+ \foned u_{xy}^2 + u_{yy}^2).
\end{align*}
Further, integrating by parts and using the fact that $f_1$ approximates $\foned \approx x$
for $x>1$, we note that
\[
\int \foned \phi_{xy}^2 \leq C \int \ftwod \phi_{xx}^2 + C \int \phi_{yy}^2.
\]
Therefore, by Gronwall's inequality
\begin{align*}
  & \sup_{0 \leq t \leq T} \int (\ftwod(\cdot,t) u_{xx}^2 + \foned(\cdot,t) u_{xy}^2 + u_{yy}^2) +
  3 \int_0^T\int \left[(\ftwod)_x u_{xxx}^2 + (\foned)_x u_{xxy}^2\right] \leq C
\end{align*}
where $C$ does not depend on $\delta$ but only on $T$ and the norm of $\phi \in X^1(\mathbb R^2)
\cap \widetilde H_x^2(W_{0\;2\;0})$.
Consequently, we can pass to the limit and conclude that
\begin{align*}
  & \sup_{0 \leq t \leq T} \int (f_2(\cdot,t) u_{xx}^2 + f_1(\cdot,t) u_{xy}^2 + u_{yy}^2) +
  3 \int_0^T\int \left[g_2 u_{xxx}^2 + g_1 u_{xxy}^2\right] \leq C.
\end{align*}

\vskip 5pt
\nd {\bf The case $j=3$.}

\vskip 5pt
 (a) {\bf The subcase $\alpha = (3,0)$.}
We choose our weight functions such that $\fthreed$ approximates
$f_3 \in W_{0\;3\;0}$.  Applying $\partial_x^3$ to (2.1), multiplying by $\fthreed u_{xxx}$
and integrating over $\mathbb R^2$, we have
\begin{align*}
& 2 \int \fthreed u_{xxx}u_{xxxt} + 2 \int \fthreed u_{xxx}u_{xxxxxx} - 2 \int \fthreed
u_{xxx}u_{xxyy} \\
& \qquad + 2 \int \fthreed u_{xxx}u_{xxxx} + 2 \int \fthreed u_{xxx}(uu_x)_{xxx} = 0.
\end{align*}
Integrating by parts gets
\begin{align*}
  & \partial_t \int \fthreed u_{xxx}^2 + 3 \int (\fthreed)_x u_{xxxx}^2 \\
  & \qquad = \int (\fthreed)_x u_{xxy}^2 + \int [\partial_t \fthreed + \partial_x^3
  \fthreed + \partial_x \fthreed]u_{xxx}^2 - 2 \int \fthreed u_{xxx}(uu_x)_{xxx} \\
  & \qquad \leq \int (\fthreed)_x u_{xxy}^2 + c \int \fthreed u_{xxx}^2 + 2 \left|
  \int \fthreed u_{xxx} (uu_x)_{xxx} \right|.
\end{align*}
Moreover
\begin{align*}
  \left|\int \fthreed u_{xxx}(uu_x)_{xxx} \right|
    & = \left|\int \fthreed [3u_{xx}^2 + 4u_x u_{xxx} + uu_{xxxx}] u_{xxx} \right| \\
  & \leq 3 \left|\int \fthreed u_{xx}^2 u_{xxx} \right| + 4 \left|\int \fthreed u_x u_{xxx}^2
  \right| + \left|\int \fthreed u u_{xxx} u_{xxxx} \right| \\
  & \leq 3 \left|\int \fthreed u_{xx}^2 u_{xxx}\right| + c ||u_x||_{L^\infty}\int \fthreed
  u_{xxx}^2 + c \left|\int (\fthreed u)_x u_{xxx}^2 \right| \\
  & \leq 3 \left|\int \fthreed u_{xx}^2 u_{xxx} \right| + c ||u||_{X^0} \int \fthreed u_{xxx}^2 \\
  & \qquad + c (||u||_{L^\infty(\mathbb R^2)} + ||u_x||_{L^\infty(\mathbb R^2)})\int
  \fthreed u_{xxx}^2 \\
  & \leq C \left|\int \fthreed u_{xx}^2 u_{xxx}\right| + c||u||_{X^0(\mathbb R^2)}\int \fthreed
  u_{xxx}^2.
\end{align*}
Now we estimate the first term on the right-hand side.

\vskip 5pt
\nd{\it Case: $x > 1$.}
Let $A_1 = \left\{x \in \mathbb R: x > 1\right\} \times \mathbb R
\subseteq \mathbb R^2$.
\begin{align*}
  \left|\int_{A_1} \fthreed u_{xx}^2 u_{xxx} \right|
    & = C \left|\int_{A_1} (\fthreed)_x u_{xx}^3 \right| \\
    & \leq C ||u_{xx}||_{L^\infty(A_1)} \int_{A_1} (\fthreed)_x u_{xx}^2 \\
    & \leq C \left(\int_{A_1} u_{xx}^2 + u_{xxxx}^2 + u_{xxy}^2\right)^{1/2}
    \int_{A_1} \ftwod u_{xx}^2 \\
    & \leq \epsilon \int_{A_1} (u_{xx}^2 + u_{xxxx}^2 + u_{xxy}^2) + C\left(\int_{A_1}
    \ftwod u_{xx}^2\right)^2.
\end{align*}
Now the terms involving $u_{xx}$ have been bounded by the previous step in the induction.
Therefore, we conclude that
\[
\left|\int_{A_1} \fthreed u_{xx}^2 u_{xxx} \right| \leq C + \epsilon \int_{A_1} (u_{xxxx}^2
+ u_{xxy}^2).
\]

\vskip 5pt
\nd {\it Case: } $x < -1$.  Let $A_{-1} = \left\{x \in \mathbb R: x < -1\right\} \times
\mathbb R \subseteq \mathbb R^2$.  We use the fact that $\fthreed \approx c$ to show
\begin{align*}
  \left|\int_{A_{-1}} \fthreed u_{xx}^2 u_{xxx} \right|
  & \leq c \left|\int_{A_{-1}} u_{xx}^2 u_{xxx} \right| \\
  & \leq c \left(\int_{A_{-1}} u_{xx}^4\right)^{1/2}\left(\int_{A_{-1}} u_{xxx}^2\right)^{1/2} \\
  & \leq C \left(\int_{A_{-1}}[u_{xx}^2 + u_{xxx}^2
  + u_{xy}^2]\right)\left(\int_{A_{-1}} u_{xxx}^2 \right)^{1/2} \\
  & \leq c ||u||_{X^0}^3.
\end{align*}
Combining these estimates for the subcase $\alpha = (3,0)$, yields
\begin{align*}
  & \int \fthreed(\cdot,t) u_{xxx}^2 + 3 \int_0^t \int (\fthreed)_x u_{xxxx}^2 \\
  & \leq
  \int \fthreed(\cdot,0) \phi_{xxx}^2 + C + \int_0^t\int (\fthreed)_x u_{xxy}^2 + c
  \int_0^t \int \fthreed u_{xxx}^2 + \epsilon \int_0^t \int_{A_1} u_{xxxx}^2 \\
  & \leq
  \int \fthreed(\cdot,0) \phi_{xxx}^2 + C + \int_0^t\int (\fthreed)_x u_{xxy}^2 + c
  \int_0^t \int \fthreed u_{xxx}^2 + \epsilon \int_0^t \int (\fthreed)_x u_{xxxx}^2
\end{align*}
Therefore,
\begin{align*}
  & \int \fthreed(\cdot,t) u_{xxx}^2 + 3 \int_0^t \int (\fthreed)_x u_{xxxx}^2 \\
    & \leq  \int \fthreed(\cdot,0) \phi_{xxx}^2 + C + \int_0^t\int (\fthreed)_x u_{xxy}^2 + c
  \int_0^t \int \fthreed u_{xxx}^2 \\
  & \leq \int \fthreed(\cdot,0) \phi_{xxx}^2 + C + C \int_0^t \int \ftwod u_{xxy}^2 +
  C \int_0^t \int \fthreed u_{xxx}^2.
\end{align*}

\vskip 5pt
 (b) {\bf The subcase $\alpha = (2,1)$.}  In this case, we take $\ftwod$ approximating
$f_2 \in W_{0\;2\;0}$.  Apply $\partial_x^2 \partial_y$ to (2.1), multiply by $\ftwod u_{xxy}$
and integrating over $\mathbb R^2$, we have
\begin{align*}
  & \partial_t \int \ftwod u_{xxy}^2 + 3 \int (\ftwod)_x u_{xxxy}^2 \\
  & \leq \int (\ftwod)_x u_{xyy}^2 + \int \ftwod u_{xxy}^2 + 2\left|\int
  \ftwod u_{xxy} (uu_x)_{xyy}\right|.
\end{align*}
Now
\begin{align*}
  \left|\int \ftwod u_{xxy}(uu_x)_{xxy} \right| & = \left|\int \ftwod
  u_{xxy} [3u_{xy}u_{xx} + 3u_x u_{xxy} + u_y u_{xxx} + uu_{xxxy}]\right| \\
  & \leq C \int \ftwod u_{xxy}u_{xy} u_{xx} + \int \ftwod u_{xxy}u_y u_{xxx}
  + C||u||_{X^0} \int \ftwod u_{xxy}^2.
\end{align*}
The second term on the right-hand side satisfies
\begin{align*}
\int \ftwod u_{xxy} u_y u_{xxx} & \leq C ||u_y||_{L^\infty} \int \ftwod (u_{xxx}^2 + u_{xxy}^2) \\
& \leq C ||u||_{X^1} \int (\fthreed u_{xxx}^2 + \ftwod u_{xxy}^2).
\end{align*}
For the first term on the right-hand side, we consider two cases.  First, for $A_1$,
\begin{align*}
  \int_{A_1} \ftwod u_{xxy}u_{xy} u_{xx} & \leq ||u_{xy}||_{L^\infty(A_1)} \left(\int_{A_1} \ftwod
  u_{xxy}^2\right)^{1/2} \left(\int_{A_1} \ftwod u_{xx}^2\right)^{1/2} \\
  & \leq \epsilon \left(\int_{A_1} u_{xy}^2 + u_{xxxy}^2 + u_{xyy}^2\right) + C
   \int_{A_1} \ftwod u_{xxy}^2 \\
   & \leq C + \epsilon \int (\ftwod)_x u_{xxxy}^2 + C \int \foned u_{xyy}^2 + C
   \int \ftwod u_{xxy}^2.
\end{align*}

Then for $A_{-1}$,
\begin{align*}
  \int_{A_{-1}} \ftwod u_{xxy} u_{xy} u_{xx} & = \int_{A_{-1}} u_{xxy} u_{xy} u_{xx} \\
  & = C \int_{A_{-1}} u_{xy}^2 u_{xxx} \\
  & \leq C \left(\int_{A_{-1}} u_{xy}^4\right)^{1/2} \left(\int_{A_{-1}} u_{xxx}^2\right)^{1/2} \\
  & \leq C \left(\int_{A_{-1}} u_{xy}^2 + u_{xxy}^2 + u_{yy}^2\right)\left(\int_{A_{-1}} u_{xxx}^2
  \right)^{1/2} \\
  & \leq C(||u||_{X^0}) \left(1 + \int \ftwod u_{xxy}^2\right)
\end{align*}

Combining these estimates and integrating with respect to $t$, we have
\begin{align*}
\int \ftwod(\cdot,t) u_{xxy}^2 + 3 \int_0^t \int (\ftwod)_x u_{xxxy}^2
& \leq C + \int \ftwod(\cdot,0) \phi_{xxy}^2 +
\epsilon \int_0^t \int (\ftwod)_x u_{xxxy}^2 \\
& \qquad + \int_0^t \int \ftwod u_{xxy}^2
+ C \int_0^t \int \foned u_{xyy}^2.
\end{align*}
Therefore,
\[
\int \ftwod(\cdot,t) u_{xxy}^2 + 3 \int_0^t \int (\ftwod)_x u_{xxxy}^2
 \leq C + \int \ftwod(\cdot,0) \phi_{xxy}^2  + \int_0^t \int \ftwod u_{xxy}^2
+ C \int_0^t \int \foned u_{xyy}^2.
\]

 (c) {\bf The subcase $\alpha = (1,2)$.}
We take our weight function $\foned$
approximating $f_1 \in W_{0\;1\;0}$.  Applying $\partial_x\partial_y^2$ to (2.1), multiplying
by $\foned u_{xyy}$ and integrating over $\mathbb R^2$, we get
\begin{align*}
  & \partial_t \int \foned u_{xyy}^2 + 3 \int (\foned)_x u_{xxyy}^2 \\
  & \leq c \int (\foned)_x u_{yyy}^2 + \int \foned u_{xyy}^2 + 2 \left|\int \foned
  u_{xyy} (uu_x)_{xyy}\right|.
\end{align*}

Now
\begin{align*}
  2 \left|\int \foned u_{xyy} (uu_x)_{xyy} \right| & = 2 \left|\int \foned u_{xyy}
  (2u_{xy}^2 + 2u_x u_{xyy} + u_{yy}u_{xx} + 2u_y u_{xxy} + uu_{xxyy})\right|.
\end{align*}
The first term on the right-hand side satisfies
\[
\int \foned u_{xyy}u_{xy}^2 = C \int \foned (u_{xy}^3)_y = 0,
\]
since $(\foned)_y = 0$.
Integrating by parts, it is clear that the second and fifth terms on the right-hand side
are bounded by
\[
C ||u||_{X^0}\int \foned u_{xyy}^2.
\]
The fourth term on the right-hand side is bounded by
\[
C||u_y||_{L^\infty} \left(\int \foned u_{xyy}^2 + \int \foned u_{xxy}^2\right).
\]
For the third-term on the right-hand side, we consider the cases when $x > 1$ and $x < -1$
separately.  First, for $x > 1$, we have
\begin{align*}
  \int_{A_1} \foned u_{xyy} u_{yy} u_{xx} & \leq ||u_{yy}||_{L^\infty(A_1)} \left(\int_{A_1}
  \foned u_{xyy}^2\right)^{1/2} \left(\int_{A_1} \foned u_{xx}^2\right)^{1/2} \\
  & \leq C \left(\int_{A_1} u_{yy}^2 + u_{xxyy}^2 + u_{yyy}^2\right)^{1/2} \left(\int_{A_1}
  \foned u_{xyy}^2\right)^{1/2} \\
  & \leq C + \epsilon \int_{A_1} u_{xxyy}^2 + C \int_{A_1} \foned u_{xyy}^2 + C \int_{A_1} u_{yyy}^2
\end{align*}
where we have used the fact that $\int \foned u_{xx}^2$ was bounded on the previous
step of the induction.  We will bound the $\epsilon$ term back on the left-hand side.
For $x < -1$, we have
\begin{align*}
  \int_{A_{-1}} \foned u_{xyy}u_{yy}u_{xx} & \approx \int \chi_{[x < -1]}u_{yy}^2 u_{xxx} \\
  & \leq \left(\int_{A_{-1}} u_{yy}^4\right)^{1/2} \left(\int_{A_{-1}} u_{xxx}^2\right)^{1/2} \\
  & \leq C\int_{A_{-1}} u_{yy}^2 + u_{xyy}^2 + (\partial_x^{-1} u_{yy})^2 \\
  & \leq C
\end{align*}
where $C$ depends only on the norm of $u \in X^1(\mathbb R^2)$.
Combining these estimates and integrating with respect to $t$, we have
\begin{align*}
  \int \foned(\cdot,t) u_{xyy}^2 + 3 \int_0^t \int (\foned)_x u_{xxyy}^2 \leq
  \int \foned(\cdot,0) \phi_{xyy}^2 + C \int_0^t \int \foned u_{xyy}^2 + C \int_0^t \int u_{yyy}^2.
\end{align*}

\vskip 5pt
 (d) {\bf The subcase $\alpha = (0,3)$.}
In this case we apply $\partial_y^3$ to (2.1), multiply by $u_{yyy}$ and integrate over
$\mathbb R^2$.  We have
\[
\partial_t \int u_{yyy}^2 \leq \int u_{yyy}^2 + 2 \left|\int u_{yyy}(uu_x)_{yyy}\right|.
\]
Now
\[
\left|\int u_{yyy}(uu_x)_{yyy} \right| = \left|\int u_{yyy}(u_xu_{yyy} + 3u_{yy}u_{xy} +
3u_y u_{xyy} + uu_{xyyy})\right|
\]
Integrating by parts as necessary, we see that the first and fourth terms on the right-hand side
are bounded by
\[
||u||_{X^0} \int u_{yyy}^2.
\]
The second term on the right-hand side is bounded by
\begin{align*}
& \left(\int u_{yy}^4\right)^{1/4} \left(\int u_{xy}^4\right)^{1/4} \left(\int u_{yyy}^2
\right)^{1/2} \\
& \leq \left(\int u_{yy}^2 + u_{xyy}^2 + (\partial_x^{-1}u_{yyy})^2\right)^{1/2}
\left(\int u_{xy}^2 + u_{xxy}^2 + u_{yy}^2\right)^{1/2} \left(\int u_{yyy}^2\right)^{1/2} \\
& \leq C||u||_{X^1(\mathbb R^2)}\left(\int u_{yyy}^2 \right)^{1/2} \\
& \leq C + C \int u_{yyy}^2.
\end{align*}
The third term on the right-hand side is bounded by
\begin{align*}
& C||u_y||_{L^\infty} \left(\int u_{xyy}^2\right)^{1/2} \left(\int u_{yyy}^2\right)^{1/2} \\
& \leq C \left(\int u_y^2 + u_{xxy}^2 + u_{yy}^2 \right)^{1/2} \left(\int u_{xyy}^2\right)^{1/2}
\left(\int u_{yyy}^2\right)^{1/2} \\
& \leq C ||u||_{X^1}^2 \left(\int u_{yyy}^2\right)^{1/2} \\
& \leq C + C \int u_{yyy}^2.
\end{align*}

Combining the estimates above and integrating with respect to $t$, we have
\begin{align*}
  \int u_{yyy}^2 \leq \int \phi_{yyy}^2 + C + C \int_0^t \int u_{yyy}^2
\end{align*}
where $C$ depends only on the norm of $\phi$ in $X^1(\mathbb R^2)$.

Combining the estimates above for (a), (b), (c) and (d) and applying Gronwall's
inequality, we conclude that
\begin{align*}
  & \sup_{0 \leq t \leq T}\int (\fthreed u_{xxx}^2 + \ftwod u_{xxy}^2 + \foned u_{xyy}^2 + u_{yyy}^2)
  \\
  & + 3 \int_0^T \int ((\fthreed)_x u_{xxxx}^2 + (\ftwod)_x u_{xxxy}^2 + (\foned)_x u_{xxyy}^2)
  \leq C
\end{align*}
where $C$ does not depend on $\delta$, but only on $T$ and the norm of $\phi \in X^1(\mathbb R^2)
\cap \widetilde H_x^3(W_{0\;3\;0})$.  Consequently, we can
pass to the limit and conclude that
\[
\sup_{0 \leq t \leq T} \int (f_3 u_{xxx}^2 + f_2 u_{xxy}^2 + f_1 u_{xyy}^2 + u_{yyy}^2) +
  3 \int_0^T \int (g_3 u_{xxxx}^2 + g_2 u_{xxxy}^2 + g_1 u_{xxyy}^2)
  \leq C.
\]

\vskip 5pt
\nd {\bf The case: $j \geq 4$.}

In this case, we take $f_{\alpha_1,\delta}$
approximating $\faone \in W_{0\;\alpha_1\;0}$.  We apply $\pa$
to (2.1), multiply by $f_{\alpha_1,\delta}\pa u$ and
integrate over $\mathbb R^2$.  We need to get a bound on
\[
\left|\int_0^T \int f_{\alpha_1,\delta} (\pa u) \pa (uu_x)\right|.
\]
In Lemma 6.2 below, we prove that
\begin{equation}
  \left|\int_0^t \int \fad (\pa u)\pa(uu_x)\right| \leq C + C \sum_{\gamma_1 + \gamma_2 =
  j} \int_0^t \int \fgd (\partial^\gamma u)^2
\end{equation}
where $C$ depends only on terms bounded in the previous step of the induction.  Consequently,
we have that
\begin{equation}
  \sup_{0\leq t\leq T} \sum_{|\alpha|=j} f_{\alpha_1,\delta} (\pa u)^2 + \sum_{|\alpha|=j} \int_0^T \int
  (f_{\alpha_1,\delta})_x (\pa u_x)^2 \leq C,
\end{equation}
where $C$ does not depend on $\delta$, but only on $T$ and the norm of $\phi \in X^1(\mathbb R^2)
\cap \widetilde H_x^j(W_{0\;j\;0})$.  Passing to the limit, we
get the desired estimate, namely,
\begin{equation}
  \sup_{0\leq t\leq T} \sum_{|\alpha|=j} \faone (\pa u)^2 + \sum_{|\alpha|=j} \int_0^T \int
  \gaone (\pa u_x)^2 \leq C.
\end{equation} \hfill $\square$

\nd{\bf Theorem 6.2}
{\it Let $f_{\alpha_1,\delta}$ approximate $\faone \in W_{0\;\alpha_1\;0}$.
Let $j = |\alpha|$, $4 \leq j \leq K$.  The following inequality holds:
\begin{equation}
\label{persist-higher}
  \left|\int_0^t \int \faoned (\pa u)\pa (uu_x)\right| \leq C + C \sum_{|\gamma| =
  j} \int_0^t \int \fgoned (\partial^\gamma u)^2
\end{equation}
for $0 \leq t \leq T$ where $C$ depends only on
\begin{align}
& \sup_{0 \leq t \leq T} \int \fgd (\partial^\gamma u)^2 \\
& \int_0^T \int (\fgd)_x (\partial^\gamma u_x)^2
\end{align}
for $\gamma = (\gamma_1,\gamma_2) \in \mathbb Z^+ \times \mathbb Z^+$, $|\gamma| \leq j-1$.
}

\vskip 5pt
\nd{\bf Proof. } In order to get bounds on the left-hand side of \eqref{persist-higher},
we use the fact that every term in the integrand is of the form
\begin{equation}
  \fad (\pa u)(\partial^r u)(\partial^s u_x)
\end{equation}
where $r_i + s_i = \alpha_i$.
Before showing the bounds on each of the terms
in the integrand we point one bound we will be using frequently:
\begin{equation}
  ||\partial^\gamma u||_{L^\infty} \leq C \left(\int (\partial^\gamma u)^2 + (\partial^\gamma u_{xx})^2
  + (\partial^\gamma u_y)^2\right)^{1/2}.
\end{equation}

\vskip 5pt
\nd {\bf The case $|s| = j$.}  In this case, $|r| = 0$ and $s = \alpha$.  Therefore,
\begin{align*}
  \int \faoned (\pa u)(\partial^r u)(\partial^s u_x) & = \int \faoned u (\pa u)(\pa u_x) \\
  & = - \frac 12 \int [\faoned u]_x (\pa u)^2 \\
  & \leq C ||u||_{X^0} \int \faoned (\pa u)^2.
\end{align*}

\vskip 5pt
\nd {\bf The case $|s| = j-1$.}  Therefore, $|r| = 1$.  We have two subcases below:

(a) {\bf The subcase $r = (1,0)$.}  In this case,
\begin{align*}
  \int \faoned (\pa u)(\partial^r u)(\partial^s u_x) & = \int \faoned (\pa u)u_x (\pa u) \\
  & \leq C||u_x||_{L^\infty} \int \faoned (\pa u)^2.
\end{align*}

(b) {\bf The subcase $r = (0,1)$.}  In this case,
\begin{align*}
  \int \faoned (\pa u)(\partial^r u)(\partial^s u_x) & = \int \faoned (\pa u)u_y
  (\partial_x^{\alpha+1} \partial_y^{\beta-2} u) \\
  & \leq C||u_y||_{L^\infty} \left(\int \faoned (\pa u)^2\right)^{1/2} \left(\int \faoned (\partial_x
  ^{\alpha_1+1}\partial_y^{\alpha_2-1}u)^2\right)^{1/2} \\
  & \leq C||u||_{X^1} \left[\int \faoned (\pa u)^2 + \int f_{\alpha+1,\delta}(\partial_x^{\alpha_1+1,
  \alpha_2-1}u)^2\right].
\end{align*}

\vskip 5pt
\nd{\bf The case $|s| = j-2$.}  We consider three subcases below:

(a) {\bf The subcase $r = (2,0)$.}  In this case
\begin{align*}
  \int \faoned (\pa u)(\partial^r u)(\partial^s u_x) & = \int \faoned (\pa u)u_{xx}
  (\partial_x^{\alpha_1-1}
  \partial_y^{\alpha_2} u) \\
  & \leq C||\foned u_{xx}||_{L^\infty} \left(\int \faoned (\pa u)^2\right)^{1/2}
  \left(\int f_{\alpha_1-1,\delta}(\partial_x^{\alpha_1-1}\partial_y^{\alpha_2}u)^2\right)^{1/2}.
\end{align*}
Now
\[
\foned u_{xx} \approx \left\{
\begin{array}{ll}
x u_{xx} & \quad x > 1 \\
u_{xx} & \quad x < -1.
\end{array}
\right.
\]
For $x < -1$, we use the fact that
\begin{align*}
||u_{xx}||_{L^\infty(A_{-1})} & \leq C \left(\int_{A_{-1}} (u_{xx}^2 + u_{xxxx}^2 + u_{xxy}^2)
\right)^{1/2} \\
& \leq C ||u||_{X^1}.
\end{align*}

For $x > 1$, we use the fact that
\begin{align*}
\foned u_{xx} & = (\foned u)_{xx} - (\foned)_{xx} u - 2(\foned)_x u_x \\
& = (\foned u)_{xx} - (\foned)_{xx} u - 2 ((\foned)_xu)_x + 2 (\foned)_{xx} u.
\end{align*}
Therefore,
\begin{align*}
||\foned u_{xx} ||_{L^\infty} & \leq ||(\foned u)_{xx}||_{L^\infty} + C ||((\foned)_x u)_x||_{L^\infty}
+ C||(\foned)_{xx}u||_{L^\infty} \\
& \leq C \left(\int ((\foned u)_{xx})^2 + ((\foned u)_{xxxx})^2 + ((\foned u)_{xxy})^2\right)^{1/2} \\
& + C \left(\int (((\foned)_xu)_x)^2 + (((\foned)_x u)_{xxx})^2 + (((\foned)_x u)_{xy})^2
\right)^{1/2} \\
& + C\left(\int ((\foned)_{xx}u)^2 + (((\foned)_{xx}u)_{xx})^2 + (((\foned)_{xx}u)_y)^2
\right)^{1/2} \\
& \leq C + C\int u^2 + u_x^2 + \foned u_{xx}^2 + u_{xxx}^2 + \foned u_{xxxx}^2 + u_y^2 + u_{xy}^2
+ \foned u_{xxy}^2 \\
& \leq C + C \sum_{|\gamma| \leq j} \int \fgd (\partial^\gamma u)^2.
\end{align*}

(b) {\bf The subcase $r=(1,1)$.}  Then
\begin{align*}
\int \faoned (\pa u)(\partial^r u)(\partial^su_x) & = \int \faoned (\pa u)u_{xy} (\partial_x^{\alpha_1}
\partial_y^{\alpha_2-1} u) \\
& \leq ||u_{xy}||_{L^\infty} \left(\int \faoned (\pa u)^2\right)^{1/2} \left(\int \faoned
(\partial_x^{\alpha_1} \partial_y^{\alpha_2-1}u)^2\right)^{1/2} \\
& \leq C||u||_{X^1} \left(\int \faoned (\pa u)^2\right)^{1/2} \left(\int \fad
(\partial_x^{\alpha_1} \partial_y^{\alpha_2-1}u)^2\right)^{1/2} \\
& \leq C + C \sum_{|\gamma| \leq j} \int \fgd (\partial^\gamma u)^2
\end{align*}
where $C$ depends only on the bounds in the statement of the theorem.

(c) {\bf The subcase $r=(0,2)$.}  Then
\begin{align*}
  \int \faoned (\pa u)(\partial^r u)(\partial^su_x) & = \int \faoned (\pa u)u_{yy}
  (\partial_x^{\alpha_1+1}\partial_y^{\alpha_2-2}u) \\
  & \leq ||u_{yy}||_{L^\infty} \left(\int \faoned (\pa u)^2\right)^{1/2} \left(
  \int \faoned (\partial_x^{\alpha_1+1} \partial_y^{\alpha_2-2} u)^2\right)^{1/2} \\
  & \leq C ||u_{yy}||_{L^\infty} \left(\int \faoned (\pa u)^2\right)^{1/2} \\
  & \leq C \left(\int u_{yy}^2 + u_{xxyy}^2 + u_{xyy}^2\right)^{1/2} \left(\int \faoned
  (\pa u)^2\right)^{1/2} \\
  & \leq C + C\sum_{|\gamma| = j}\int \fgd (\partial^\gamma u)^2
\end{align*}
where $C$ depends only on the bounds in the statement of the theorem.

\vskip 5pt
\nd {\bf The case $|s| = j-3$ for $j \geq 5$.}

In this case we consider $x > 1$ and $x < -1$ separately.
First, for $x < -1$, we have
\begin{align*}
  \int_{A_{-1}} \faoned (\pa u)(\partial^r u)(\partial^s u_x)
  & \leq \int \left(\int_{A_{-1}} (\pa u)^2
  \right)^{1/2} \left(\int_{A_{-1}} (\partial^r u)^4\right)^{1/4}
  \left(\int_{A_{-1}} (\partial^s u_x)^4\right)^{1/4}.
\end{align*}
Now
\[
\left(\int_{A_{-1}} (\partial^r u)^4\right)^{1/4} \leq \left(\int_{A_{-1}} (\partial^r u)^2
+ (\partial^r u_x)^2 + (\partial^r u_y)^2\right)^{1/2}
\]
and $|r| = 3$ implies each of these terms is bounded by $C$ where $C$ depends only on
\[
\int \fgd (\partial^\gamma u)^2 \qquad \text{for } |\gamma| \leq j-1
\]
since $j \geq 5$ and each of these terms has derivatives of order $\leq 4$.
Also,
\[
\left(\int_{A_{-1}} (\partial^s u_x)^4 \right)^{1/4} \leq \left(\int_{A_{-1}}
(\partial^su_x)^2 + (\partial^su_{xx})^2 + (\partial^s u_{xy})^2 \right)^{1/2}
\]
and $|s| = j-3$.  Therefore, each of these terms has order at most $j-1$ and thus
bounded by
\[
C \sum_{|\gamma| \leq j-1} \int \fgd (\partial^\gamma u)^2.
\]
Therefore, for $x < -1$, we have
\[
\int_{A_{-1}} \faoned (\pa u)(\partial^r u)(\partial^s u_x)
\leq C + C \int_{A_{-1}} \faoned (\pa u)^2,
\]
where $C$ depends only on the terms in the statement of the theorem.

Now for $x > 1$, we have
\begin{align*}
\int_{A_1} \faoned (\pa u)(\partial^r u)(\partial^s u_x)
& \leq \left(\int_{A_1} \faoned (\pa u)^2
\right)^{1/2} \left(\int_{A_1} (f_{r_1/2,\delta} \partial^ru)^4\right)^{1/4} \\
& \times \left(\int_{A_1} (f_{(s_1+1)/2,\delta}
\partial^su_x)^4\right)^{1/4}.
\end{align*}
since $r_1 + s_1 = \alpha_1$.
Now
\begin{align*}
\left(\int_{A_1} (f_{r_1/2,\delta} \partial^r u)^4 \right)^{1/4} & \leq
\left(\int_{A_1} (f_{r_1/2,\delta} \partial^r u)^2 + ((f_{r_1/2,\delta} \partial^ru)_x)^2
+ ((f_{r_1/2,\delta} \partial^r u)_y)^2\right)^{1/2} \\
& \leq C \left(\int_{A_1} f_{r_1,\delta} (\partial^ru)^2 + f_{r_1,\delta} (\partial^r u_x)^2
+ f_{r_1,\delta}(\partial^ru_y)^2\right)^{1/2}.
\end{align*}
Since $|r| = 3$, each of the terms above has order at most four, and, therefore,
is bounded by a constant $C$ which depends only on
\[
\int \fgd (\partial^\gamma u)^2 \qquad \text{for } |\gamma| \leq 4 \leq j-1,
\]
since here we are assuming $j \geq 5$.

Similarly,
\begin{align*}
  & \left(\int (f_{(s_1+1)/2,\delta}(\partial^s u_x))^4\right)^{1/4} \\
  & \leq C
  \left(\int (f_{(s_1+1)/2,\delta}(\partial^s u_x))^2+((f_{(s_1+1)/2,\delta}(\partial^su_x))_x)^2
  + ((f_{(s_1+1)/2,\delta}(\partial^su_x))_y)^2\right)^{1/2} \\
  & \leq C \left(\int f_{s_1+1,\delta}[(\partial^su_x)^2 + (\partial^su_{xx})^2+ (\partial^s
  u_{xy})^2]\right)^{1/2}.
\end{align*}
Since $|s| = j-3$ each of these terms is of order at most $j-1$.  Therefore, each of
these terms is bounded by
\[
C \sum_{|\gamma| \leq j-1} \int \fgd (\partial^\gamma u)^2.
\]

\vskip 5pt
\nd {\bf The case $|s| = j-3$ when $j = 4$.}

In this case, $|s| = 1$.  Therefore,
either $s = (1,0)$ or $s = (0,1)$.  For $s = (1,0)$, we have
\begin{align*}
\int \faoned (\pa u)(\partial^r u) (\partial^s u_x) & = \int \faoned (\pa u)(\partial^r u)
u_{xx} \\
& \leq \left(\int \faoned (\pa u)^2\right)^{1/2} \left(\int (f_{r_1/2,\delta} \partial^r u)^4
\right)^{1/4} \left(\int (\foned u_{xx})^4\right)^{1/4}.
\end{align*}

Now
\begin{align*}
  \left(\int (f_{r_1/2,\delta}(\partial^r u))^4\right)^{1/4} & \leq \left(\int f_{r_1,\delta}
  [(\partial^r u)^2+ (\partial^r u_x)^2 + (\partial^r u_y)^2]\right)^{1/2}.
\end{align*}
We note that $|r| = 3$.  Therefore, each of these terms is at most of order $4 = j$.
Further,
\begin{align*}
  \left(\int (\foned u_{xx})^4\right)^{1/4} & \leq C \left(\int \ftwod[u_{xx}^2 + u_{xxx}^2
  + u_{xxy}^2]\right)^{1/2}.
\end{align*}
We note that each of these terms is of order at most $j-1$.
Combining these estimates, we conclude that
\[
\int \faoned (\pa u)(\partial^r u) (\partial^s u_x) \leq C + C \sum_{|\gamma| = j}
\int \fgd (\partial^\gamma u)^2,
\]
where $C$ depends only on
\[
\sum_{|\gamma| \leq j-1} \int \fgd (\partial^\gamma u)^2.
\]

The case in which $s = (0,1)$ is handled similarly.

\vskip 5pt
\nd{\bf The case $|s| \leq j-4$.}

In this case, we bound the terms as follows:
\begin{align*}
\int \faoned (\pa u)(\partial^r u)(\partial^s u_x) & \leq ||f_{(s_1+1)/2,\delta}\partial^su_x
||_{L^\infty} \left(\int \faoned(\pa u)^2\right)^{1/2} \left(\int f_{r_1,\delta} (\partial^r u)^2
\right)^{1/2}.
\end{align*}
Now
\begin{align*}
  ||f_{(s_1+1)/2,\delta}(\partial^su_x)||_{L^\infty} & \leq C\left(\int f_{s_1+1,\delta}[
  (\partial^s u_x)^2 + (\partial^s u_{xxx})^2 + (\partial^s u_{xy})^2]\right)^{1/2}.
\end{align*}
Since $|s| \leq j-4$, all of these terms are of order at most $j-1$.  Therefore,
each of these terms is bounded by
\[
C \sum_{|\gamma| \leq j-1} \int \fgd (\partial^\gamma u)^2
\]
and, therefore,
\[
\int \faoned (\pa u) (\partial^r u)(\partial^s u_x) \leq C + C \sum_{|\gamma| = j}
\int \fgd (\partial^\gamma u)^2.
\]
where $C$ depends only on
\[
C \sum_{|\gamma| \leq j-1} \int \fgd (\partial^\gamma u)^2.
\]\hfill $\square$

\renewcommand{\theequation}{\thesection.\arabic{equation}}
\setcounter{equation}{0}\section{Main Theorem} In this section we
state and prove our main theorem, which states that if the initial
data $\phi$ possesses certain regularity
and sufficient decay at infinity, then the solution $u(t)$
will be smoother than $\phi$.  In particular if the initial data satisfies
\[
\int \phi^2 + (1 + x_+^L) (\partial_x^L \phi)^2 + (\partial_y^L \phi)^2 < \infty,
\]
then the solution will {\it gain} $L$ derivatives in $x$.  More specifically,
\[
\int_0^T \int t^{L-1}(1+e^{\sigma x_-}) (\partial_x^{2L} u)^2 < \infty
\]
for $\sigma > 0$ arbitrary, where $T$ is the existence time of the solution.
\\



\nd{\bf Theorem 7.1} (Main Theorem). {\it Let $T>0$ and let $u$ be the
solution of \eqref{e201} in the region
$[0,\,T]\times\mathbb{R}^{2}$ such that $u\in
L^{\infty}([0,\,T]:\,{\cal Z}_{L})$ for some $L\geq 2.$ Then}
\begin{eqnarray}
\label{e701}\sup_{0\leq t\leq
T}\int_{\mathbb{R}^{2}}f_\alpha\,(\pa u)^{2}\,dx\,dy +
\int_{0}^{T}\int_{\mathbb{R}^{2}}g_\alpha\,(\pa u_x)^{2}\,dx\,dy\,dt \leq C
\end{eqnarray}
{\it for $L+1\leq |\alpha| \leq 2L - 1$, $2L-|\alpha|-\alpha_2 \geq 1$ where
$f_\alpha \in
W_{\sigma,\;2L - |\alpha|-\alpha_2,\;|\alpha|-L}$ and $g_\alpha\in
W_{\sigma,\;2L - |\alpha| - \alpha_2 - 1,\;|\alpha|-L},$ $\sigma>0$
arbitrary.}\\
\\
{\it Proof.} By assumption, $u\in L^{\infty}([0,\,T]:\,{\cal
Z}_{L}).$ Therefore $u_{t}\in
L^{\infty}([0,\,T]:\,L^{2}(\mathbb{R}^{2})),$ then $u\in
C([0,\,T]:\,L^{2}(\mathbb{R}^{2}))\cap C_{w}([0,\,T]:\,{\cal
Z}_{L}).$ Hence $u:[0,\,T]\rightarrow {\cal Z}_{L}$ is a weakly
continuous function. In particular, $u(\,\cdot\,,\,\cdot\,,t)\in
{\cal Z}_{L}$ for every $t.$ Let $t_{0}\in (0,\,T)$ and
$u(\,\cdot\,,\,\cdot\,,t_{0})\in {\cal Z}_{L},$ then there are
$\{\phi^{(n)}\}\subseteq C_{0}^{\infty}(\mathbb{R}^{2})$ such that
$\partial_{x}^{-1}\phi_{yy}^{(n)}$ are in
$C_{0}^{\infty}(\mathbb{R}^{2})$ and
$\phi^{(n)}(\,\cdot\,,\,\cdot\,)\rightarrow
u(\,\cdot\,,\,\cdot\,,t_{0})$ in ${\cal Z}_L.$ Let $u^{(n)}$ be the
unique solution of \eqref{e201} with initial data
$\phi^{(n)}(x,\,y)$ at time $t=t_{0}.$ By Corollary 4.4, the
solution $u^{(n)}\in L^{\infty}([t_{0},\,t_{0} +
\delta]:\,X^{1}(\mathbb{R}^{2}))$ for a time interval $\delta$ which
not depend on $n.$ By Theorem 6.1, $u^{(n)}\in
L^{\infty}([t_{0},\,t_{0} + \delta]:\,{\cal Z}_L)$ and
\begin{eqnarray}
\label{e702}\int_{t_{0}}^{t_{0} +
\delta}\int_{\mathbb{R}^{2}}g_{\alpha_1}\,(\pa
u_x)^{2}\,dx\,dy\,dt\leq C
\end{eqnarray}
for $|\alpha|=L$, $\alpha_1 \neq 0$,
where $g_{\alpha_1}\in W_{\sigma,\;\alpha_1 - 1,\;0}$ and $C$ depends only on the norm
of $\phi^{(n)}\in {\cal Z}_{L}.$ Also by Theorem 6.1, we have
(non-uniform) bounds on
\begin{eqnarray}
\label{e703}\sup_{[t_{0},\,t_{0} + \delta]}\sup_{(x,y)}\left[(1 +
x^{+})^{k}|\partial_x^{\alpha_1} u^{(n)}(x,\,y,\,t)| + |\partial_y^{\alpha_2}
u^{(n)}(x,\,y,\,t)|\right] <+\infty
\end{eqnarray}
for each $n,\,k$ and $\alpha_1,\,\alpha_2.$ Therefore, the a priori
estimates in Lemma 5.1, are justified for each $u^{(n)}$ in the
interval $[t_{0},\,t_{0} + \delta].$\\
\\
We start our induction with $|\alpha|=L+1$.  In this case, we take $g_\alpha
\in W_{\sigma,\;L - 2 - \alpha_2,\;1}$ and let $f_\alpha =\frac{1}{3}\int_{-\infty}^{x}
g_\alpha(z,\,t)\,dz.$  We note that $2L-|\alpha|-\alpha_2 \geq 1$ by assumption.
Therefore, $L-2-\alpha_2 \geq 0$.  As
shown in Lemma 5.1, we have the following bounds on the higher
derivatives of $u^{(n)},$
\begin{eqnarray}
\label{e704}\sup_{[t_{0},\,t_{0} +
\delta]}\int_{\mathbb{R}^{2}}f_\alpha \,(\pa u^{(n)})^{2}\,dx\,dy
+ \int_{t_{0}}^{t_{0} +
\delta}\int_{\mathbb{R}^{2}}g_\alpha \,(\pa u^{(n)}_x)^{2}\,dx\,dy\,dt\leq
C
\end{eqnarray}
where $C$ depends only on the norm of $u^{(n)}\in
L^{\infty}([0,\,T]:\,{\cal Z}_{L})$ and the term in \eqref{e702}. We
conclude, therefore, that the constant $C$ in \eqref{e704} depends
only on $||\phi^{(n)}||_{{\cal Z}_{L}}.$ We continue this procedure
inductively. For the $|\alpha|^{th}$ step, let $g_\alpha \in W_{\sigma,\;2L - |\alpha|
-\alpha_2 - 1,\;|\alpha|-L}$ for $\alpha_2 \leq 2L-|\alpha|-1$ and define
$f_\alpha =\frac{1}{3}\int_{-\infty}^{x}g_\alpha(z,\,t)\,dz.$ The non-uniform
bounds on $u^{(n)}$ in \eqref{e702} allows us to use Lemma 5.1 and
our inductive hypothesis to conclude that
\begin{eqnarray*}
\sup_{[t_{0},\,t_{0} +
\delta]}\int_{\mathbb{R}^{2}}f_\alpha \,(\pa u^{(n)})^{2}\,dx\,dy + \int_{t_{0}}^{t_{0} +
\delta}\int_{\mathbb{R}^{2}}g_\alpha \,(\pa u^{(n)}_x)^{2}\,dx\,dy\,dt\leq C
\end{eqnarray*}
where again $C$ does not depend on $n,$ but only on the norm of
$\phi^{(n)}\in {\cal Z}_{L}.$ By Corollary 4.5,
\begin{eqnarray*}
u^{(n)}\stackrel{*}\rightharpoonup u\quad \mbox{weakly in}\quad
L^{\infty}([t_{0},\,t_{0} + \delta]:\,X^{1}(\mathbb{R}^{2})).
\end{eqnarray*}
Therefore, we can pass to the limit and conclude that
\begin{eqnarray}
\label{e706}\sup_{[t_{0},\,t_{0} +
\delta]}\int_{\mathbb{R}^{2}}f_\alpha \,(\pa u)^{2}\,dx\,dy + \int_{t_{0}}^{t_{0} +
\delta}\int_{\mathbb{R}^{2}}g_\alpha \,(\pa u)^{2}\,dx\,dy\,dt\leq C.
\end{eqnarray}
This proof is continued inductively up to $|\alpha| =2L - 1.$ Since $\delta$
does not depend on $n,$ this result is valid over the whole interval
$[0,\,T].$ \hfill $\square$

\end{document}